\documentclass[a4paper]{article}

\usepackage[english]{babel}
\usepackage{amsfonts}
\usepackage[T1]{fontenc}

\usepackage[a4paper,top=3cm,bottom=2cm,left=3cm,right=3cm,marginparwidth=1.75cm]{geometry}

\usepackage{amsmath}
\usepackage{amsthm}
\usepackage{graphicx}
\usepackage{subcaption}
\usepackage{tikz}
\usepackage[utf8]{inputenc}
\usepackage{pgfplots} 
\usepackage{pgfgantt}
\usepackage{pdflscape}
\pgfplotsset{compat=newest} 
\pgfplotsset{plot coordinates/math parser=false} 
\newlength\fwidth
\newlength\fheight
\newcommand{\D}{\ensuremath{\mathbf{D}}}
\newcommand{\Q}{\ensuremath{\mathsf{Q}}}
\newcommand{\T}{\ensuremath{\mathsf{T}}}
\newcommand{\Qbar}{\ensuremath{\bar{\mathsf{Q}}}}
\newcommand{\Real}{\mathbb{R}}
\newcommand{\bzero}{\mathbf{0}}

\newcommand{\da}{\Delta a}

\newcommand{\q}{\mathbf q}
\newcommand{\p}{\mathbf p}
\newcommand{\qbar}{\bar{\mathbf q}}
\newcommand{\pbar}{\bar{\mathbf p}}
\newcommand{\qdot}{\dot{\mathbf {q}}}
\newtheorem{theorem}{Theorem}
\usepackage[colorinlistoftodos]{todonotes}
\usepackage[colorlinks=true, allcolors=blue]{hyperref}
\theoremstyle{definition}
\newtheorem{remark}{Remark}

\title{Performance Assessment of Energy-preserving,
Adaptive Time-step Variational Integrators}

\author{Harsh Sharma \thanks{Corresponding Author, Email Address for Correspondence: hasharma@ucsd.edu} \thanks{Postdoctoral Scholar, Department of Mechanical and Aerospace Engineering, University of California San Diego}, Jeff Borggaard \thanks{Professor, Department of Mathematics, Virginia Tech}, Mayuresh Patil \thanks{Professor of Practice, Department of Aerospace Engineering, Georgia Tech}, and Craig Woolsey \thanks{Professor, Kevin T. Crofton Department of Aerospace and Ocean Engineering, Virginia Tech}}
\begin{document}
\maketitle
\begin{abstract} 
A fixed time-step variational integrator cannot preserve momentum, energy, and symplectic form simultaneously for nonintegrable systems. This barrier can be overcome by treating time as a discrete dynamic variable and deriving adaptive time-step variational integrators that conserve the energy in addition to being symplectic and momentum-preserving. Their utility, however, is still an open question due to the numerical difficulties associated with solving the discrete governing equations. In this work, we investigate the numerical performance of energy-preserving, adaptive time-step variational integrators. First, we compare the time adaptation and energy performance of the energy-preserving adaptive algorithm with the adaptive variational integrator for Kepler's two-body problem. Second, we apply tools from Lagrangian backward error analysis to investigate numerical stability of the energy-preserving adaptive algorithm. Finally, we consider a simple mechanical system example to illustrate the backward stability of this energy-preserving, adaptive time-step variational integrator. 
\end{abstract}
%
%
%
%
\section{Introduction} 
Variational integrators are a class of structure-preserving integrators, that are derived by discretizing the action principle rather than the governing differential equation.  The basic idea behind these methods is to obtain an approximation of the action integral called the discrete action. Stationary points of the discrete action give discrete time trajectories of the mechanical system. These integrators have good long-time energy behavior and they conserve the invariants of the dynamics in the presence of symmetries.  The basics of variational integrators can be reviewed in  \cite{Lew2016AIntegrators, sharma2020review}  and more detailed theory can be found in \cite{Marsden2001DiscreteIntegrators,wendlandt1997mechanical}. A more general framework encompassing variational integrators, asynchronous variational integrators, and symplectic-energy-momentum integrators is discussed in \cite{leok2005generalized}. Due to their symplectic nature,  variational integrators are ideal for long-time simulation of conservative or weakly dissipative systems \cite{kane2000variational} found in astrophysics and molecular dynamics. The discrete trajectories obtained using variational integrators display excellent energy behavior for exponentially long times. \par
The fixed time-step variational integrators derived from the discrete variational principle cannot preserve the energy of the system exactly.  To conserve the energy in addition to preserving the symplectic structure and conserving the momentum, time adaptation needs to be used.  These energy-preserving, adaptive time-step variational integrators were first developed for conservative systems in \cite{Kane1999Symplectic-energy-momentumIntegrators}  by imposing an additional energy preservation equation to compute the time-step. The same integrators were derived through a variational approach for a more general case of time-dependent Lagrangian systems in \cite{Marsden2001DiscreteIntegrators} by preserving the discrete energy obtained from the discrete variational principle in the extended phase space. Guibout and Bloch \cite{guibout2004discrete} have discussed the subtle differences between the underlying discrete variational principles. Recently, Sharma et al. \cite{sharma2018energy,sharma2019acc,sharma2019cdc} extended these methods to systems with external nonconservative forcing and rigid body motion. 
\par
Despite their excellent conservation properties, the utility of these methods is still an open question due to the difficulties associated with their numerical implementation. These adaptive algorithms are implicit and require solving a coupled nonlinear system of equations at every time-step to update both the configuration variables and the time variable. Unlike traditional adaptive algorithms, the adaptive time-step computation in these methods is inherently coupled to the discrete dynamics. In fact, existence of solutions for these discrete trajectories is not always guaranteed. Shibberu \cite{Shibberu2006Is,shibberu2005regularize} has discussed the well-posedness of these adaptive algorithms and identified points in the extended state space where the adaptive algorithm has no solutions. Even for points where solutions exist, the governing update equations are ill-conditioned. This introduces a discrete energy error at every iteration and can accumulate
in long-time simulations. Therefore, the backward stability of energy-preserving, adaptive time-step variational integrator  algorithms needs to be investigated.
\par 
The early development of backward error analysis in 1960s was motivated by numerical linear algebra applications and it was mainly used as a tool for evaluating the propagation of rounding errors in matrix algorithms. Although Wilkinson \cite{wilkinson1960error} first used it in the context of eigenvalue problems, the underlying idea of backward error analysis has also been applied to numerical integration of dynamical systems \cite{eirola1993aspects,hairer1997life}. Truncation errors introduced by the numerical integrator play the same role as rounding errors in the numerical linear algebra setting. The focus in the backward error analysis is to answer the question: what dynamical system does the numerical integrator solve? Thus, instead of studying the difference between the trajectories obtained using the numerical integrator and the exact solution, we look for a nearby dynamical system which would be solved exactly by the numerical integrator. For dynamical systems with underlying geometric structure, we would like the modified system to also possess similar features. The existence of such a nearby dynamical system has direct implications for the accurate long-time behavior of the simulation.
\par
Sanz-Serna \cite{sanz1996backward} first applied tools from backward error analysis to symplectic integrators and showed that numerical trajectories from symplectic integrators can be interpreted as the exact solution of a perturbed Hamiltonian system. Reich \cite{reich1999backward} used a recursive definition of modified vector fields to provide a unifying framework for the backward error analysis of geometric numerical integration (GNI) methods. Hairer \cite{hairer1997variable} developed variable time-step symplectic integrators and proved backward stability of the adaptive algorithm by using time transformations. Recently, Vermeeran \cite{vermeeren2017modified} studied the backward stability of fixed time-step variational integrators from the Lagrangian perspective. 
\par
%
%
%
The main goal of this work is to understand the numerical performance of energy-preserving, adaptive time-step variational integrators. First, we provide a comparison between different adaptive time-stepping approaches through a numerical example of Kepler's two-body problem. We also obtain a conservative bound on the discrete energy error that is accumulated from inexact solutions to the residual of the discrete energy equation. This shows that the error can be controlled through the residual tolerance and the numerical precision that is used. Variable precision arithmetic in the energy-preserving adaptive algorithm is used to show how the energy performance can be improved when using more significant digits in the computation. Motivated by the improvement in energy performance, we theoretically justify these results by investigating the stability of the adaptive algorithm. We first use a time-transformation to represent the adaptive algorithm as a fixed time-step algorithm and then apply backward error analysis tools in the transformed setting.\par
%
%
The remainder of the paper is organized as follows. In Section \ref{s:2}, we review the basics of extended Lagrangian mechanics and derive the energy-preserving, adaptive time-step variational integrators. In Section \ref{s:3}, we compare the numerical performance of energy-preserving, adaptive time-step variational integrators and adaptive variational integrators for the Kepler's two-body problem, with special focus on their energy behavior and time adaptivity. We also study the energy behavior of the adaptive algorithm and the effect of higher precision arithmetic on the conservation properties.  In Section \ref{s:6}, we use tools from backward error analysis to study algorithmic stability of energy-preserving, adaptive time-step variational integrators. We first introduce a time transformation to represent the energy-preserving, adaptive time-step variational integrator as a fixed time-step variational integrator. We then use these backward stability concepts to derive modified equations and the corresponding modified Lagrangian for a simple mechanical system. Finally, in Section \ref{s:7} we provide concluding remarks and suggest future research directions.
%
%
%
%
\section{Background} 
 \label{s:2} 
In this section, we review the basics of extended Lagrangian mechanics and the derivation of energy-preserving, adaptive time-step variational integrators. Drawing on the work of Marsden et al.~\cite{Kane1999Symplectic-energy-momentumIntegrators, Marsden2001VariationalMechanics, Marsden2001DiscreteIntegrators}, we first derive equations of motion in continuous-time from the variational principle and then derive the corresponding variational integrators by considering the discretized variational principle in the discrete-time domain.  To achieve this, we first derive extended Euler-Lagrange equations from Hamilton's principle. After deriving equations of motion, we use concepts of extended discrete mechanics developed in \cite{Marsden2001DiscreteIntegrators} to derive extended discrete Euler-Lagrange equations and write these in a time-marching form to obtain energy-preserving, adaptive time-step variational integrators for Lagrangian systems. 
\subsection{Extended Lagrangian Mechanics}
\label{sc:21}
Hamilton's principle of stationary action is one of the most fundamental results of classical mechanics and is commonly used to derive equations of motion for a variety of systems. The forces and interactions that govern the dynamical evolution of the system are easily determined through Hamilton's principle in a formulaic and elegant manner.  Hamilton's principle \cite{Goldstein1980ClassicalMechanics} states that: The motion of the system between two fixed points from $t_0$ to $t_f$ is such that the action integral has a stationary value for the actual path of the motion. To derive the Euler-Lagrange equations via Hamilton's principle, we start by defining the configuration space, tangent space and path space.
\par
Consider a time-dependent Lagrangian system with configuration manifold $\Q$ and time space $\Real$. In the extended Lagrangian mechanics framework \cite{Marsden2001DiscreteIntegrators}, we treat time as a dynamic variable and define the extended configuration manifold $\Qbar=\Real \times \Q$; the corresponding state space $\T\Qbar$ is $\Real \times \T\Q$, where the tangent bundle $\T\Q$ is the union of all tangent spaces to $\Q$. The extended Lagrangian is $L: \Real \times \T\Q \to \Real$. \par
In the extended Lagrangian mechanics framework, $t$ and $\q$ are both parametrized by an independent variable $a$. The two components of a trajectory $c$ are $c(a)=(c_t(a),c_{\q}(a))$. The extended path space is 
\begin{equation}
\bar{\mathcal{C}}=\left\{c:[a_0,a_f] \to \Qbar \ |\  \mbox{$c$ is a $C^2$ curve and $c'_t(a) >0$} \right\}.
\end{equation}
For a given path $c(a)$, the initial time is $t_0=c_t(a_0)$ and the final time is $t_f=c_t(a_f)$. The extended action $\bar{\mathfrak{B}} : \bar{\mathcal{C}} \to \Real$ is
\begin{equation}
\bar{\mathfrak{B}} = \int ^{t_f} _{t_0} L(t, \q(t), \qdot(t)) dt.
\label{eq:action}
\end{equation}
Since time is a dynamic variable in this framework, we substitute $(t,\q(t),\dot{\q}(t))=\left(c_t(a),c_{\q}(a),\frac{c'_{\q}(a)}{c'_t(a)} \right)$ in the above equation to get 
\begin{equation}
\bar{\mathfrak{B}} = \int ^{a_f} _{a_0} L\left(c_t(a),c_{\q}(a),\frac{c'_{\q}(a)}{c'_t(a)} \right)c'_t(a) \  da.
\end{equation}
We compute the first variation of the action as
\begin{equation}
\delta \bar{\mathfrak{B}} = \int^{a_f}_{a_0} \left[ \frac{\partial L}{\partial t} \delta c_t + \frac{\partial L}{\partial \q}\cdot \delta c_{\q}  +  \frac{\partial L}{\partial \qdot}\cdot \left( \frac{\delta c'_{\q}(a)}{c'_t}  -\frac{c'_{\q} \delta c'_t(a) }{(c'_t)^2} \right) \right] c'_t(a)da + \int ^{a_f} _{a_0} L \delta c'_t(a) \  da .
\end{equation}
Using integration-by-parts and setting the variations at the endpoints to zero gives
\begin{equation}
\delta \bar{\mathfrak{B}} = \int^{a_f}_{a_0} \left[ \frac{\partial L}{\partial \q} c'_t - \frac{d}{da} \frac{\partial L}{\partial \qdot} \right] \cdot \delta c_{\q}(a)da + \int ^{a_f} _{a_0}  \left[ \frac{\partial L}{\partial t}c'_t + \frac{d}{da}\left( \frac{\partial L}{ \partial \qdot}\cdot \frac{c'_{\q}}{c'_t} - L \right) \right] \delta c_t(a) \  da .
\end{equation}
Using $dt=c'_t(a) da$ in the above expression gives two equations of motion. The first is the Euler-Lagrange equation of motion 
\begin{equation}
\frac{\partial L}{\partial \q}   - \frac{d}{dt} \left( \frac{\partial L}{\partial \qdot}  \right) =\bzero, 
\label{eq:contp}
\end{equation}
which is the same as the equation obtained using the classical Lagrangian mechanics framework. The second equation is  
\begin{equation}
\frac{\partial L}{ \partial t} + \frac{d}{dt} \left( \frac{\partial L}{ \partial \qdot} \cdot \dot{\q} - L \right) =0,
\label{eq:conte}
\end{equation}
which describes how the energy of the system evolves with time. 
\begin{remark}
It should be noted that both \eqref{eq:contp} and \eqref{eq:conte} depend only on the associated curve $\q(t)$. The time component $c_t(a)$ of the extended path cannot be determined from the governing equations. This comes from the fact that the energy evolution equation \eqref{eq:conte} is redundant as it is a consequence of the Euler-Lagrange equation \eqref{eq:contp}. Thus, the "velocity" of time, i.e. $c'_t(a)$, is indeterminate in the continuous-time formulation.
\end{remark} 
\subsection{Energy-preserving, Adaptive Time-step Variational Integrators}
\label{sc:22}  
For the extended Lagrangian mechanics, we define the extended discrete state space $\Qbar \times \Qbar$. The {\em extended discrete path space} is 
\begin{equation}
\bar{\mathcal{C}}_d = \{\ c: \{\ 0,...,N \} \to \Qbar\ | \ \mbox{$c_t(k+1) > c_t(k)$  for all $k$}  \}\ .
\end{equation}
The {\em extended discrete action map} $\bar{\mathcal{B}}_d : \bar{\mathcal{C}}_d \to \Real$ is
\begin{equation}
\bar{\mathcal{B}}_d = \sum _{k=0}^{N-1} L_d(t_k,\q_k,t_{k+1},\q_{k+1}),
\end{equation}
where $L_d : \Qbar \times \Qbar \to \Real$ is the {\em extended discrete Lagrangian function} that approximates the action integral between two successive configurations. Taking variations of the extended discrete action map leads to
\begin{multline}
      \delta \bar{B}_d = \sum_{k=1}^{N-1} \left[ \D_4L_d(t_{k-1},\q_{k-1},t_{k},\q_{k}) + \D_2L_d(t_k,\q_k,t_{k+1},\q_{k+1}) \right] \cdot \delta \q_k  \\  +\sum_{k=1}^{N-1} \left[ \D_3L_d(t_{k-1},\q_{k-1},t_{k},\q_{k}) + \D_1L_d(t_k,\q_k,t_{k+1},\q_{k+1}) \right]  \delta t_k  = 0,
    \end{multline}
where $\D_i$ denotes differentiation with respect to the $i^{th}$ argument of the discrete Lagrangian $L_d$. Applying Hamilton's principle of least action and setting variations at endpoints to zero results in the extended discrete Euler-Lagrange equations 
\begin{equation}
 \D_4L_d(t_{k-1},\q_{k-1},t_{k},\q_{k}) + \D_2L_d(t_k,\q_k,t_{k+1},\q_{k+1}) = \bzero,
\label{eq:vi_pvv}
\end{equation}

\begin{equation}
\D_3L_d(t_{k-1},\q_{k-1},t_{k},\q_{k}) + \D_1L_d(t_k,\q_k,t_{k+1},\q_{k+1}) = 0.
\label{eq:vi_ev}
\end{equation}
Given $(t_{k-1},\q_{k-1},t_{k},\q_{k})$, the extended discrete Euler-Lagrange equations can be solved to obtain $\q_{k+1}$ and $t_{k+1}$. This  extended discrete Lagrangian system can be seen as a numerical integrator of the continuous-time nonautonomous Lagrangrian system with adaptive time-steps.
\par
In the extended discrete mechanics framework, we define the discrete momentum $\p_k$ by
\begin{equation}
\p_{k}= \D_4L_d(t_{k-1},\q_{k-1},t_{k},\q_{k}).
\end{equation}
We also introduce the discrete energy 
\begin{equation}
E_{k}=\D_3L_d(t_{k-1},\q_{k-1},t_{k},\q_{k}).
\label{eq:dE}
\end{equation}
Using the discrete momentum and discrete energy definitions, we can re-write the extended discrete Euler-Lagrange equations \eqref{eq:vi_pvv} and \eqref{eq:vi_ev} in the following form
\begin{align}
\label{eq:pimp}
-\D_2L_d(t_k,\q_k,t_{k+1},\q_{k+1}) &= \p_k,\\
\label{eq:eimp}
\D_1L_d(t_k,\q_k,t_{k+1},\q_{k+1}) &= E_k,\\
\label{eq:pexp}
\p_{k+1}&= \D_4L_d(t_k,\q_k,t_{k+1},\q_{k+1}),\\
\label{eq:eexp}
E_{k+1}&= -\D_3L_d(t_k,\q_k,t_{k+1},\q_{k+1}).
\end{align}
Given $(t_k,\q_k,\p_k,E_k)$, the coupled nonlinear equations \eqref{eq:pimp} and \eqref{eq:eimp} are solved implicitly to obtain $\q_{k+1}$ and $t_{k+1}$. The configuration $\q_{k+1}$ and time $t_{k+1}$ are then used in \eqref{eq:pexp} and \eqref{eq:eexp} to obtain $(\p_{k+1},E_{k+1})$ explicitly. The extended discrete Euler-Lagrange equations were first written in the  time-marching form in \cite{Marsden2001DiscreteIntegrators} and are also known as symplectic-energy-momentum integrators. 
%
\section{Time Transformation and Adaptivity}
\label{s:3}
In this section, we compare the numerical performance of the energy-preserving, adaptive time-step variational integrators discussed in Section \ref{sc:22} to the adaptive variational integrators from \cite{schmitt2017adaptive}. Due to their variational derivations, both classes of integrators are symplectic and momentum-preserving but motivation for the time adaptation is very different. Adaptive variational integrators are designed for computational efficiency whereas the energy-preserving, adaptive time-step variational integrators are focused on preserving the discrete energy in addition to being symplectic and momentum-preserving. 
\subsection{Adaptive Variational Integrators}
For a Hamiltonian system $H(\q, \p)$, Schmitt and Leok  \cite{schmitt2017adaptive} prescribed an {\em a priori} time transformation $t \to a$, given by $\frac{dt}{da}=g(\q, \p)$, with auxilary states $(q^t,p^t)$ to form a new Hamiltonian 
\begin{equation}
\bar{H}(\qbar,\bar{\p})=g(\q,\p) (H(\q,\p) + p^t),
\label{eq: newH}
\end{equation}
where $(\qbar,\pbar)=\left( \begin{bmatrix}
\q \\ q^t
\end{bmatrix},\begin{bmatrix}
\p \\ p^t
\end{bmatrix}\right)$. Subsequently, the physical time $t$ is treated as a configuration variable by choosing $q^t=t$ with the corresponding momentum variable $p^t=-H(\q_0,\p_0)$.  Hamilton's equations for the new Hamiltonian \eqref{eq: newH} are given by 
\begin{equation}
\frac{d\qbar}{da} =
\begin{bmatrix}
g(\q,\p) \frac{\partial H}{\partial \p} \\ g(\q,\p)
\end{bmatrix}, \quad \quad \quad  \frac{d\pbar}{da} =
\begin{bmatrix}
-g(\q,\p) \frac{\partial H}{\partial \q} \\ 0
\end{bmatrix}.
\end{equation}
Applying Hamiltonian variational integrators with fixed time-step $\Delta a$ to the transformed Hamiltonian system result in adaptive variational integrators for the original Hamiltonian system $H(\q,\p)$. We apply the implicit midpoint rule for time-integration of the transformed Hamiltonian system resulting in
\begin{equation*}
\begin{bmatrix}
\frac{\q_{k+1}-\q_k}{\da} \\
\frac{q^t_{k+1}-q^t_k}{\da}
\end{bmatrix} = \begin{bmatrix}
 g(\q_{av})H_{\p}(\q_{av},\p_{av}) \\
g(\q_{av})
\end{bmatrix}, \quad \begin{bmatrix}
\frac{\p_{k+1}-\p_k}{\da} \\
\frac{p^t_{k+1}-p^t_k}{\da}
\end{bmatrix} = \begin{bmatrix}
- g(\q_{av})H_{\q}(\q_{av},\p_{av}) \\
0
\end{bmatrix},
\end{equation*}
where $\q_{av}=\frac{\q_{k+1} + \q_k}{2}$, $\p_{av}=\frac{\p_{k+1} + \p_k}{2}$, and  we have assumed that the monitor functions only depend on the configuration, i.e. $g(\q,\p)=g(\q)$. Given $(\q_k, q^t_k, \p_k)$, we need to solve the following coupled nonlinear equations to  obtain $(\q_{k+1}, q^t_{k+1}, \p_{k+1})$
\begin{equation}
\label{eq:adap_qk}
\frac{\q_{k+1}-\q_k}{\da} =  g\left( \frac{\q_{k+1} + \q_k}{2}\right)\frac{\partial H}{\partial \p}\left(\frac{\q_{k+1} + \q_k}{2},\frac{\p_{k+1} + \p_k}{2}\right),
\end{equation}
\begin{equation}
\label{eq:adap_pk}
\frac{\p_{k+1}-\p_k}{\da} =  -g\left( \frac{\q_{k+1} + \q_k}{2}\right)\frac{\partial H}{\partial \q}\left(\frac{\q_{k+1} + \q_k}{2},\frac{\p_{k+1} + \p_k}{2}\right),
\end{equation}
\begin{equation}
\label{eq:adap_tk}
\frac{q^t_{k+1}-q^t_k}{\da} =  g\left( \frac{\q_{k+1} + \q_k}{2}\right).
\end{equation}
\subsection{Numerical Comparison for Kepler's Two-body Problem}
\label{s:32}
For the comparison of the above adaptive variational integrator to the energy-preserving, adaptive time-step variational integrator presented in Section~\ref{sc:22}, we consider Kepler's two-body problem with two different initial conditions. The Hamiltonian for this conservative system is
\begin{equation}
H(q^1,p^1,q^2,p^2)=\frac{(p^1)^2 +(p^2)^2}{2} - \frac{1}{\sqrt{(q^1)^2 +(q^2)^2}}.
\end{equation}
For initial conditions, we  consider $q^1(0)=1-e, q^2(0)=0, p^1(0)=0, p^2(0)=\sqrt{\frac{1+e}{1-e}}$ for $e=0.1$ and $0.7$, where $e$ is the eccentricity of the elliptical trajectory.
\subsubsection{Implementation Details} 
\begin{itemize}
\item For energy-preserving, adaptive time-step variational integrator (EpAVI), we use the midpoint rule to obtain the extended discrete Lagrangian
\begin{equation}
L_d=(t_{k+1}-t_k)L\left( \frac{\q_k + \q_{k+1}}{2},\frac{\q_{k+1}-\q_k}{t_{k+1}-t_k}\right). 
\label{eq:Ld_comp}
\end{equation} 
The time-marching equations for the EpAVI algorithm can be derived by substituting the discrete Lagrangian from \eqref{eq:Ld_comp} in \eqref{eq:pimp}-\eqref{eq:eexp}.
\item For adaptive variational integrator (AVI), we consider two monitor functions used in \cite{schmitt2017adaptive}. The first is the arclength parametrization
\begin{equation}
g_1(\q)= (2(H_0-V(\q)) + \nabla V(\q)^T M^{-1} \nabla V(\q) )^{-\frac{1}{2}} ,
\end{equation}
 and the second is the problem-specific monitor function 
 \begin{equation}
 g_2(\q)=\q^{\top}\q,
 \end{equation}
 which is motivated by Kepler's second law. We denote the adaptive algorithms based on the first monitor function $g_1$ and second monitor function $g_2$ by AVI1 and AVI2, respectively.
 \item All numerical studies are done with MATLAB 2020b. To compute the trajectory error plots, we compute the reference solution using the in-built ode45 solver with $reltol=1e-12$ and $abstol=1e-14$. 
 \item For a fair comparison between the two adaptive approaches, $\Delta a$ for the adaptive algorithm \eqref{eq:adap_qk}-\eqref{eq:adap_tk} is chosen such that the resulting initial time step $h_0$ is same as the one used for the energy-preserving adaptive algorithm. 
 \end{itemize}
 \subsubsection{Results}
\begin{figure}
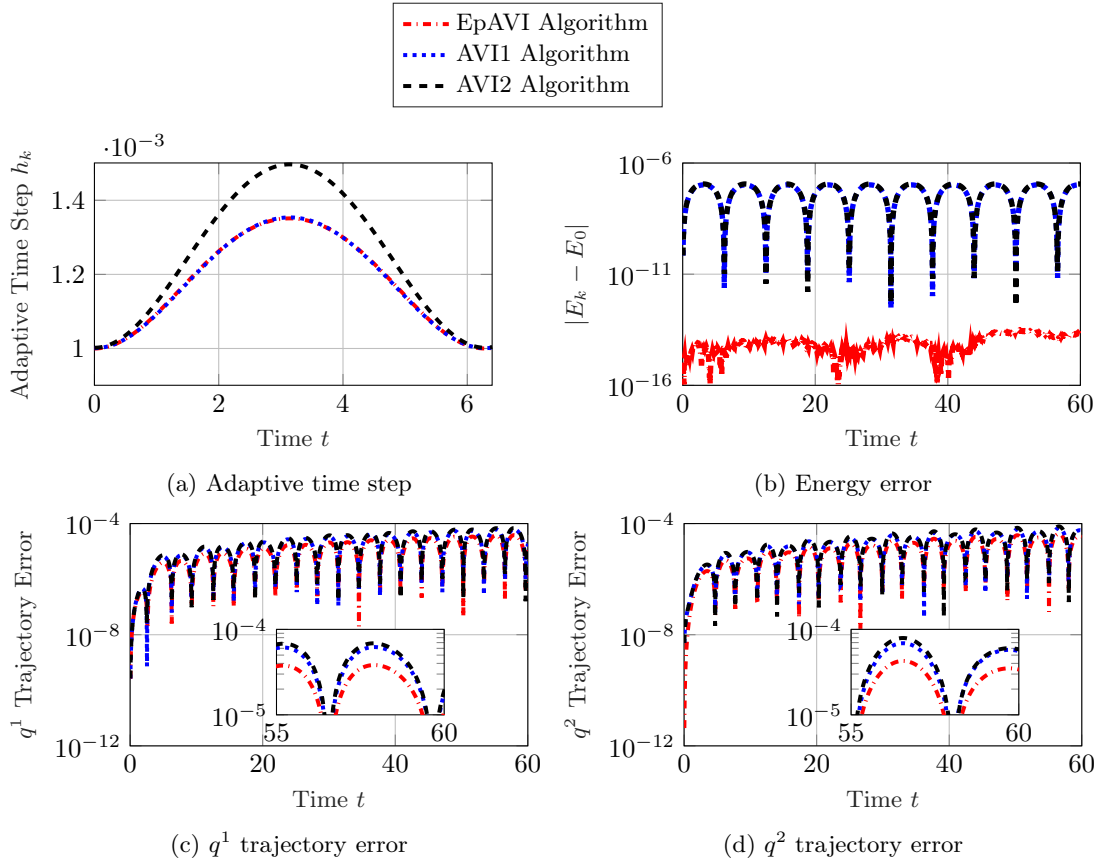

\captionsetup[subfigure]{oneside,margin={2.0cm,0 cm}}
\begin{subfigure}{.37\textwidth}
       \setlength\fheight{5.5 cm}
        \setlength\fwidth{\textwidth}
\input{Figures/tikz/H_e_01.tex}
\caption{Adaptive time step}
\label{fig: compe01h}
    \end{subfigure}
    \hspace{1.5cm}
    \begin{subfigure}{.37\textwidth}
           \setlength\fheight{5.5 cm}
           \setlength\fwidth{\textwidth}
\raisebox{-60mm}{\input{Figures/tikz/E_e_01.tex}}
\caption{Energy error}
\label{fig: compe01e}
    \end{subfigure}\\
    \begin{subfigure}{.37\textwidth}
       \setlength\fheight{5.5 cm}
        \setlength\fwidth{\textwidth}
\input{Figures/tikz/q1_e_01_inset.tex}
\caption{$q^1$ trajectory error}
\label{fig: q101h}
    \end{subfigure}
    \hspace{1.5cm}
    \begin{subfigure}{.37\textwidth}
           \setlength\fheight{5.5 cm}
           \setlength\fwidth{\textwidth}
\input{Figures/tikz/q2_e_01_inset.tex}
\caption{$q^2$ trajectory error}
\label{fig: q201e}
    \end{subfigure}
    \caption{Comparison between adaptive variational integrators for $e=0.1$}
\label{fig:compe01}
\end{figure}

We compare the energy error, trajectory error, and adaptive time-step behavior for $e=0.1$ and $e=0.7$ in Figure \ref{fig:compe01} and Figure \ref{fig:compe07}, respectively. For $e=0.1$, we observe that EpAVI algorithm and AVI1 algorithm give nearly same adaptive time-step behavior, while AVI2 algorithm achieves marginally higher adaptive time-step values. The energy error plots in Figure  \ref{fig: compe01e} demonstrate that AVI1 and AVI2 exhibit bounded energy error magnitude around $10^{-7}$ whereas EpAVI algorithm has energy error magnitude around $10^{-14}$. Trajectory error plots for $q_1$ and $q_2$ demonstrate the improved accuracy achieved by EpAVI algorithm. All three algorithms presented exhibit nearly  the same accuracy, with the EpAVI algorithm performing marginally better. 
\par 
For the higher eccentricity value $e=0.7$, the adaptive variational integrators (AVI1 and AVI2) exhibit significantly higher adaptive time-step values compared with the EpAVI algorithm. Unlike the $e=0.1$ case where the adaptive time step $h_k$ does not increase substantially from the initial time step $h_0=0.001$, all three adaptive algorithms achieve significant speedup with the adaptive time-step value $h_k$ increasing by almost 12-15 times the initial value $h_0$. Both adaptive variational integrators have energy error around $10^{-4}$ whereas the EpAVI algorithm has energy error around $10^{-13}$. Similar to the $e=0.1$ case, the AVI1 and AVI2 algorithms exhibit similar trajectory error, and the EpAVI algorithm performs marginally better.  
\begin{figure}
\captionsetup[subfigure]{oneside,margin={2.0cm,0 cm}}
\begin{subfigure}{.37\textwidth}
       \setlength\fheight{5.5 cm}
        \setlength\fwidth{\textwidth}
\input{Figures/tikz/H_e_07.tex}
\caption{Adaptive time step}
\label{fig: compe07h}
    \end{subfigure}
    \hspace{1.5cm}
    \begin{subfigure}{.37\textwidth}
           \setlength\fheight{5.5 cm}
           \setlength\fwidth{\textwidth}
\raisebox{-60mm}{\input{Figures/tikz/E_e_07.tex}}
\caption{Energy error}
\label{fig: compe07e}
    \end{subfigure} \\
    \begin{subfigure}{.37\textwidth}
       \setlength\fheight{5.5 cm}
        \setlength\fwidth{\textwidth}
\input{Figures/tikz/q1_e_07_inset.tex}
\caption{$q^1$ trajectory error}
\label{fig: q107h}
    \end{subfigure}
    \hspace{1.5cm}
    \begin{subfigure}{.37\textwidth}
           \setlength\fheight{5.5 cm}
           \setlength\fwidth{\textwidth}
\input{Figures/tikz/q2_e_07_inset.tex}
\caption{$q^2$ trajectory error}
\label{fig: q207e}
    \end{subfigure}
    \caption{Comparison between adaptive variational integrators for $e=0.7$}
\label{fig:compe07}
\end{figure}
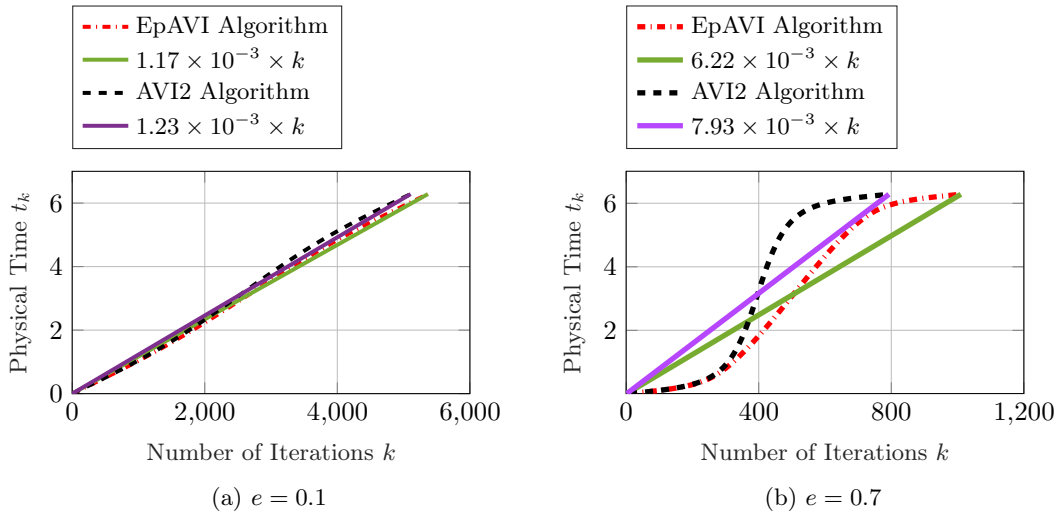
\begin{figure}[h]
\captionsetup[subfigure]{oneside,margin={1.5cm,0 cm}}
\hspace{0.25cm}
\begin{subfigure}{.37\textwidth}
       \setlength\fheight{5.5 cm}
        \setlength\fwidth{\textwidth}
%
%
\definecolor{mycolor1}{rgb}{0.46667,0.67451,0.18824}%
\definecolor{mycolor2}{rgb}{0.49412,0.18431,0.55686}%
\definecolor{mycolor3}{rgb}{1.00000,0.00000,1.00000}%
\begin{tikzpicture}

\begin{axis}[%
width=0.951\fheight,
height=0.536\fheight,
at={(0\fheight,0\fheight)},
scale only axis,
xmin=0,
xmax=6000,
xlabel style={font=\color{white!15!black}},
xlabel={\small Number of Iterations $k$},
ymin=0,
ymax=7,
ylabel style={font=\color{white!15!black}},
ylabel={\small Physical Time $t_k$},
axis background/.style={fill=white},
xmajorgrids,
ymajorgrids,
legend style={at={(0,1.1)}, anchor=south west, legend cell align=left, align=left, draw=white!15!black, font=\small}
]
\addplot [color=red, dashdotted, line width=1.5pt]
  table[row sep=crcr]{%
1	0\\
284	0.283811476570918\\
436	0.437933090229308\\
564	0.569323100648035\\
679	0.688972443994317\\
786	0.801908427532908\\
887	0.910119294010656\\
984	1.01565605609721\\
1078	1.11954528157548\\
1170	1.22284973551268\\
1260	1.32552715267047\\
1349	1.42868094616188\\
1438	1.53346765727929\\
1527	1.63989379458417\\
1616	1.74794718192334\\
1706	1.85883750999983\\
1798	1.97382450740588\\
1893	2.09421161693263\\
1991	2.22004393382849\\
2094	2.35393847476644\\
2205	2.49989515267043\\
2327	2.66198074570366\\
2471	2.85499230219193\\
2681	3.13830473331927\\
2948	3.49824115306092\\
3091	3.68932647195197\\
3212	3.84940847088546\\
3322	3.99331549808721\\
3424	4.12514430295596\\
3522	4.25017260877121\\
3616	4.36847196398412\\
3708	4.4826214343102\\
3798	4.5926622292036\\
3887	4.69985646092209\\
3976	4.80541463841928\\
4065	4.90933454327478\\
4155	5.01277265524186\\
4246	5.11570943738207\\
4338	5.21814554015964\\
4433	5.32227826952294\\
4531	5.42805389698242\\
4633	5.53650850111899\\
4741	5.64970058032759\\
4858	5.77066350587484\\
4988	5.90339067715377\\
5142	6.05891487922599\\
5366	6.28331706314657\\
};
\addlegendentry{EpAVI Algorithm}

\addplot [color=mycolor1, line width=1.5pt]
  table[row sep=crcr]{%
1	0\\
5366	6.28331706314657\\
};
\addlegendentry{$1.17 \times 10^{-3} \times k$}

\addplot [color=black, dashed, line width=1.5pt]
  table[row sep=crcr]{%
1	0\\
261	0.261335628072629\\
400	0.402725497734536\\
517	0.523342756107922\\
621	0.632157487292716\\
717	0.734201088285772\\
808	0.832554202348547\\
894	0.927122372042504\\
977	1.02002144943708\\
1057	1.1111921275824\\
1135	1.20171899585603\\
1211	1.29155634374274\\
1286	1.38185446506577\\
1360	1.47259815645884\\
1433	1.56375715837748\\
1505	1.65528657883715\\
1577	1.74843225559744\\
1650	1.84451520630682\\
1723	1.94223033326398\\
1797	2.04290876493815\\
1873	2.14794744181381\\
1951	2.2573869141861\\
2033	2.37409798696353\\
2119	2.49815223327278\\
2212	2.63395163564746\\
2317	2.78893873071229\\
2448	2.98402553103824\\
2830	3.55396810247203\\
2934	3.70686581675727\\
3027	3.84197449785552\\
3113	3.96528009473877\\
3194	4.07978958063904\\
3272	4.18842682087234\\
3348	4.29263661138157\\
3422	4.39247365452593\\
3495	4.48933532363844\\
3568	4.58455074654921\\
3640	4.67683549410231\\
3712	4.76750660547714\\
3785	4.85780547539434\\
3859	4.94769464865294\\
3934	5.03715114895749\\
4010	5.12616682405769\\
4088	5.21588867028731\\
4168	5.30627918532537\\
4251	5.39842241695169\\
4338	5.49335153493394\\
4429	5.59099700359184\\
4526	5.69342960554604\\
4631	5.80265352394235\\
4748	5.92269343188127\\
4886	6.06258718916342\\
5082	6.25947170580821\\
5106	6.28352103831457\\
};
\addlegendentry{AVI2 Algorithm}

\addplot [color=mycolor2, line width=1.5pt]
  table[row sep=crcr]{%
1	0\\
5106	6.28352103831457\\
};
\addlegendentry{$1.23 \times 10^{-3} \times k$}

\end{axis}

\begin{axis}[%
width=1.227\fheight,
height=0.658\fheight,
at={(-0.16\fheight,-0.072\fheight)},
scale only axis,
xmin=0,
xmax=1,
ymin=0,
ymax=1,
axis line style={draw=none},
ticks=none,
axis x line*=bottom,
axis y line*=left
]
\end{axis}
\end{tikzpicture}%
\caption{$e=0.1$}
\label{fig: transform1}
    \end{subfigure}
    \hspace{1.5cm}
    \begin{subfigure}{.37\textwidth}
           \setlength\fheight{5.5 cm}
           \setlength\fwidth{\textwidth}
%
%
\definecolor{mycolor1}{rgb}{0.46667,0.67451,0.18824}%
\definecolor{mycolor2}{rgb}{0.71765,0.27451,1.00000}%
\definecolor{mycolor3}{rgb}{1.00000,0.07451,0.65098}%
\begin{tikzpicture}

\begin{axis}[%
width=0.951\fheight,
height=0.536\fheight,
at={(0\fheight,0\fheight)},
scale only axis,
xmin=0,
xmax=1200,
xlabel style={font=\color{white!15!black}},
xlabel={\small Number of Iterations $k$},
xtick={0,400,800,1200},
ymin=0.00834326579261024,
ymax=7.00834326579261,
ylabel style={font=\color{white!15!black}},
ylabel={\small Physical Time $t_k$},
axis background/.style={fill=white},
xmajorgrids,
ymajorgrids,
legend style={at={(0,1.1)}, anchor=south west, legend cell align=left, align=left, draw=white!15!black, font=\small}
]
\addplot [color=red, dashdotted, line width=2.0pt]
  table[row sep=crcr]{%
9	0.0080041156071502\\
48	0.0478621351820721\\
72	0.0740422356086583\\
91	0.0963563982040796\\
107	0.116664800945614\\
122	0.137337727267322\\
135	0.156849665675963\\
147	0.176465131028067\\
158	0.196065691461399\\
168	0.2154722733884\\
177	0.234445145723498\\
186	0.255068155509889\\
194	0.274986604738046\\
202	0.296608495667101\\
209	0.31710499057067\\
216	0.339253542746292\\
222	0.359696122193327\\
228	0.381620170916221\\
234	0.405161612119628\\
240	0.430461014173716\\
246	0.457659705008155\\
251	0.481875025035038\\
256	0.507577737405768\\
261	0.534832455425999\\
266	0.563691786536765\\
271	0.594193279364958\\
276	0.626357036271088\\
281	0.660184398851129\\
286	0.695657987366189\\
291	0.732743154193486\\
296	0.771390655162236\\
301	0.811540131731931\\
306	0.853123896506759\\
312	0.904817059895322\\
318	0.958348601669059\\
324	1.01359864001302\\
330	1.0704550680307\\
336	1.12881522005853\\
343	1.19867917553006\\
350	1.27033196550974\\
357	1.34365661613288\\
364	1.41854843248757\\
371	1.49491249341463\\
379	1.58387664875613\\
387	1.67452830282639\\
395	1.76675277343816\\
404	1.8722483770083\\
413	1.97944296754781\\
422	2.0881867947943\\
432	2.21064955988595\\
443	2.34711280761258\\
455	2.49776742853123\\
468	2.66267774289156\\
484	2.86741430161851\\
509	3.18929073032155\\
534	3.51067670751695\\
549	3.70190397387285\\
562	3.8659891407799\\
573	4.00328233763867\\
583	4.12661880276539\\
593	4.24833927565714\\
602	4.35633778004842\\
611	4.46271342774742\\
620	4.56731539397674\\
628	4.65868119927438\\
636	4.74841126495915\\
644	4.83638743912252\\
651	4.91182867688997\\
658	4.98573849994045\\
665	5.05801606926593\\
672	5.12854920680445\\
678	5.18752259919677\\
684	5.24503311700755\\
690	5.30098228271038\\
696	5.35526190759424\\
702	5.40775511778315\\
708	5.45833902519018\\
713	5.49894409989406\\
718	5.53806773502163\\
723	5.57564687546437\\
728	5.61162840625002\\
733	5.64597278511258\\
738	5.67865701816061\\
743	5.70967652642332\\
748	5.73904563263545\\
753	5.76679664843107\\
758	5.79297778204545\\
763	5.81765024644824\\
768	5.8408849979453\\
774	5.8669782517062\\
780	5.89125301024399\\
786	5.91385072822186\\
792	5.93490969533673\\
799	5.9577094005399\\
806	5.97878964911445\\
813	5.99832910077907\\
821	6.01897944544748\\
829	6.03804182904435\\
838	6.05782214578949\\
848	6.07800324584969\\
858	6.09655474868532\\
869	6.11535342048205\\
881	6.13423301555099\\
894	6.15308369356092\\
909	6.17313821273262\\
926	6.19411735100812\\
945	6.21588651314312\\
968	6.24056375654072\\
1001	6.27415248751686\\
1010	6.28315754058622\\
};
\addlegendentry{EpAVI Algorithm}

\addplot [color=mycolor1, line width=2.0pt]
  table[row sep=crcr]{%
2	0.00622711351888938\\
1010	6.28315754058622\\
};
\addlegendentry{$6.22 \times 10^{-3} \times k$}

\addplot [color=black, dashed, line width=2.0pt]
  table[row sep=crcr]{%
9	0.00800670372859713\\
47	0.0468703642989112\\
70	0.0719803068694773\\
89	0.0943309616950501\\
105	0.114730255216955\\
119	0.134115187482053\\
132	0.153716343182055\\
144	0.17347824972444\\
155	0.193283558790426\\
165	0.212950335089545\\
174	0.232231986816146\\
183	0.25325044300439\\
191	0.273609278070353\\
198	0.29289683789591\\
205	0.313736777022882\\
212	0.336322188489703\\
218	0.357236254710415\\
224	0.379750050184043\\
230	0.404038440427144\\
235	0.425774841248312\\
240	0.449004488777405\\
245	0.473865004120398\\
250	0.500507272388859\\
255	0.529096532298013\\
259	0.553490646390742\\
263	0.579346219931949\\
267	0.606770400610571\\
271	0.635877706941415\\
275	0.66679026575855\\
279	0.699637979144882\\
283	0.734558597080536\\
286	0.76219667195312\\
289	0.791146692058533\\
292	0.821475198723306\\
295	0.853250925500902\\
298	0.886544624793714\\
301	0.921428841860688\\
304	0.957977629298284\\
307	0.996266194852069\\
310	1.03637047532061\\
313	1.07836662944089\\
316	1.12233044302377\\
319	1.16833664031901\\
322	1.21645809668678\\
325	1.26676494920321\\
328	1.31932360389987\\
331	1.37419564096911\\
334	1.43143662250327\\
337	1.49109481118114\\
339	1.53222978690212\\
341	1.57446619406153\\
343	1.6178111600392\\
345	1.66226943961999\\
347	1.70784322274676\\
349	1.7545319446258\\
351	1.80233210045139\\
353	1.85123706719412\\
356	1.9266434226887\\
359	2.00446437998244\\
362	2.08463337584078\\
365	2.16706459452359\\
368	2.25165243163303\\
371	2.33827127484494\\
374	2.4267756408226\\
377	2.5170006990968\\
380	2.60876320262309\\
384	2.73315631058551\\
388	2.8594085215334\\
393	3.01898929900119\\
406	3.43514784612933\\
410	3.56124179643962\\
414	3.68542414151375\\
417	3.77699670865729\\
420	3.86700687447478\\
423	3.95527373173809\\
426	4.04163493490898\\
429	4.12594757686861\\
432	4.20808866083348\\
435	4.28795518875347\\
438	4.36546389823684\\
441	4.44055068821785\\
443	4.48923970684007\\
445	4.53682236204224\\
447	4.58329265581926\\
449	4.62864767412589\\
451	4.6728874070792\\
453	4.71601456258907\\
455	4.75803437583158\\
457	4.7989544167865\\
460	4.85829418993217\\
463	4.91522260250781\\
466	4.96978894698816\\
469	5.02204964312205\\
472	5.07206696763012\\
475	5.11990787429454\\
478	5.16564291231771\\
481	5.20934524705956\\
484	5.25108978410935\\
487	5.29095239508672\\
490	5.3290092415715\\
493	5.36533619207091\\
496	5.40000832588828\\
499	5.43309951709762\\
502	5.46468209147974\\
505	5.49482654918825\\
508	5.52360134601838\\
511	5.55107272641658\\
515	5.58578381416021\\
519	5.61843584854194\\
523	5.64916581014404\\
527	5.67810283385188\\
531	5.70536829191531\\
535	5.73107596877651\\
539	5.75533230435985\\
544	5.78376254095986\\
549	5.8102592895234\\
554	5.83498662345107\\
559	5.85809434838882\\
564	5.87971909062378\\
570	5.90388575907957\\
576	5.92628987188436\\
582	5.9471050141168\\
589	5.96958709448984\\
596	5.99033519460647\\
603	6.00954101163961\\
611	6.02981723778771\\
619	6.04851809134846\\
628	6.06790936792493\\
638	6.08768105214313\\
648	6.10584696159538\\
659	6.12424830174859\\
671	6.1427257631268\\
685	6.1625376259766\\
700	6.18207378637806\\
717	6.20255460220358\\
737	6.2249611705472\\
762	6.2512456028129\\
793	6.28250441744001\\
};
\addlegendentry{AVI2 Algorithm}

\addplot [color=mycolor2, line width=2.0pt]
  table[row sep=crcr]{%
2	0.00793245507247775\\
793	6.28250441744001\\
};
\addlegendentry{$7.93 \times 10^{-3} \times k$}

\end{axis}

\begin{axis}[%
width=1.227\fheight,
height=0.658\fheight,
at={(-0.16\fheight,-0.072\fheight)},
scale only axis,
xmin=0,
xmax=1,
ymin=0,
ymax=1,
axis line style={draw=none},
ticks=none,
axis x line*=bottom,
axis y line*=left
]
\end{axis}
\end{tikzpicture}%
\caption{$e=0.7$}
\label{fig: transform2}
    \end{subfigure}
    \caption{Time transformation $t(a)$ for $h_0=0.001$}
    \label{fig: transform}
\end{figure}
\par 
Based on the numerically computed discrete time values $\{t_k\}_{k=0}^{N}$,  we have compared the time transformation $t(a)$ in Figure \ref{fig: transform}  for both adaptive approaches. These plots reveal the numerical relation between the physical time $t$ and the transformed time $a$ over one cycle for both eccentricity cases. In Figure \ref{fig: transform1}, both energy-preserving adaption and monitor function approaches have a similar $t(a)$ relation for $e=0.1$. For an initial time-step of $h_0=0.001$, the AVI2 algorithm achieves an average time-step of $1.23h_0$ whereas the average time-step for the EpAVI algorithm is $1.17h_0$. By contrast, Figure \ref{fig: transform2} shows how both monitor function and energy-preserving adaptation approaches have significantly higher average time-step for $e=0.7$. The AVI2 algorithm achieves an average time-step of $7.93h_0$ whereas the average time-step for the EpAVI algorithm is $6.22h_0$. Furthermore, $t(a)$ plots in Figure \ref{fig: transform} numerically validate the bounded nature of $t'(a)$, i.e. $ 0<t'(a)\leq const$. The fact that $t'(a)$ remains bounded will be used as an assumption in the backward error analysis in Section \ref{s:6}.
\begin{figure}
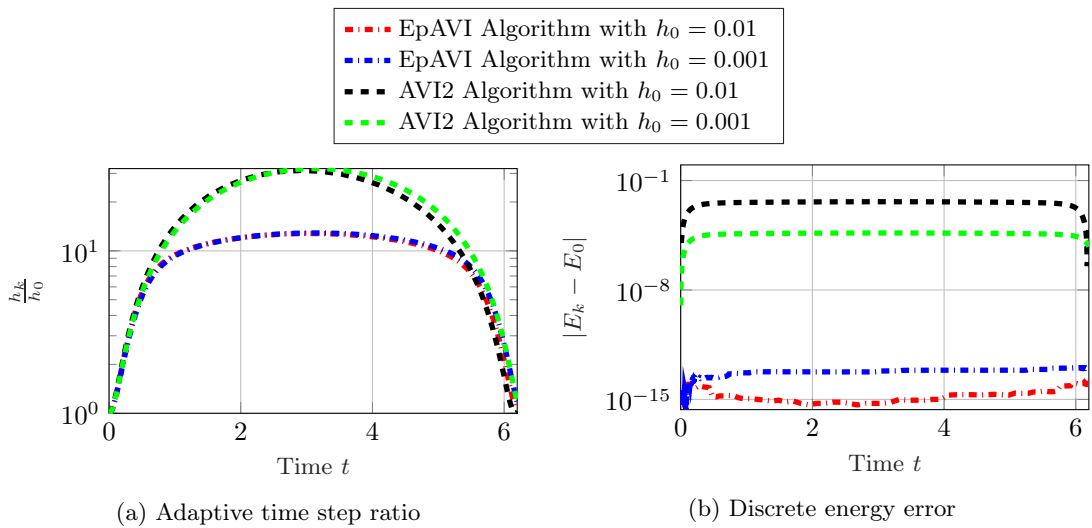

\captionsetup[subfigure]{oneside,margin={1.5cm,0 cm}}
\begin{subfigure}{.37\textwidth}
       \setlength\fheight{5.5 cm}
        \setlength\fwidth{\textwidth}
\input{Figures/tikz/Hrat_e_07.tex}
\caption{Adaptive time step ratio}
\label{fig: rat07h}
    \end{subfigure}
    \hspace{1.5cm}
    \begin{subfigure}{.37\textwidth}
           \setlength\fheight{5.5 cm}
           \setlength\fwidth{\textwidth}
\raisebox{-62mm}{\input{Figures/tikz/Erat_e_07.tex}}
\caption{Discrete energy error}
\label{fig: rat07e}
    \end{subfigure}
    \caption{Effect of initial time-step $h_0$ on adaptive variational integrators for $e = 0.7$}
\label{fig:rat07}
\end{figure}
\par 
We have also studied the effect of initial time-step $h_0$ on the numerical performance of the adaptive algorithms. For the AVI2 algorithm, increasing the initial time-step from $h_0=0.001$ to $h_0=0.01$ leads to an increase in the energy error. The EpAVI algorithm, on the other hand, exhibits small decrease in the energy error for the smaller initial time-step $h_0=0.001$. This unexpected decrease in energy error for larger time-steps is due to residual error in not solving the discrete energy equation exactly. We discuss the dependence of the residual error on energy behavior of the EpAVI algorithm in detail in Section \ref{s:4}. Both algorithms demonstrate similar numerical behavior of $\frac{h_k}{h_0}$ for both  $h_0$ values which shows that the speedup achieved up the adaptive algorithms is independent of initial time-step and is mainly governed by the discrete dynamics. 
\begin{remark}
We observe slow linear error growth for all three adaptive algorithms in the trajectory error plots for both eccentricity values in Figure \ref{fig:compe01} and Figure \ref{fig:compe07}. This error growth is due to the integrable nature of Kepler's two-body problem. Completely integrable systems form an interesting class of Hamiltonian systems, and symplectic integrators are known to exhibit linear error growth for such systems. More details about this result can be found in \cite{hairer2006geometric}. \end{remark}
\subsection{Energy Behavior}
\label{s:33}
As shown in Section \ref{sc:22}, the extended discrete equations for the adaptive time-step variational integrators conserve the discrete energy $E_k$ exactly. Unfortunately, numerical algorithms do not solve the discrete energy equation exactly and the residual from approximating solutions to this equation affects the errors in subsequent iterations. It is possible that the accumulation of the discrete energy error leads to a significant energy error that overcomes the energy-preservation advantage of the adaptive time-stepping algorithm. Thus, we need to investigate the long-term behavior of energy error.
\par
For a given initial configuration $(t_0,\q_0,\p_0,E_0)$, the discrete energy $E_k$ after $k$ adaptive time-steps can be written as the telescoping series
\begin{equation}
E_k= E_0 + \sum_{i=1}^k  (E_i -E_{i-1}).
\end{equation}
We can use this expression to obtain a conservative bound on the discrete energy error 
\begin{equation}
| E_k- E_0| \leq  \sum_{i=1}^k | E_k -E_{k-1}| \leq k \max_i |E_i -E_{i-1} |,
\label{eq: dEbound}
\end{equation}
where $\max_i|E_i-E_{i-1}|$ is the maximum discrete energy error introduced in one iteration during the entire numerical simulation. Using this conservative bound, we can compare the energy behavior of fixed time-step variational integrators with the worst case scenario for the EpAVI algorithm. For example, if we consider an example where the fixed time-step algorithm exhibits discrete energy error around $10^{-6}$ then the EpAVI algorithm, with $10^{-15}$ accuracy in solving discrete energy equation, will exhibit better energy performance for $10^{9}$ steps
\begin{equation}
    k=\frac{10^{-6}}{\max_i |E_i -E_{i-1} |}= \frac{10^{-6}}{10^{-15}}=10^{9}  \ steps .
\end{equation}
\begin{figure}[h]
\captionsetup[subfigure]{oneside,margin={1.5cm,0 cm}}
\begin{subfigure}{.37\textwidth}
       \setlength\fheight{5.5 cm}
        \setlength\fwidth{\textwidth}
%
%
\begin{tikzpicture}

\begin{axis}[%
width=0.951\fheight,
height=0.536\fheight,
at={(0\fheight,0\fheight)},
scale only axis,
xmin=0,
xmax=6.19836391244326,
xlabel style={font=\color{white!15!black}},
xlabel={\small Time $t$},
ymin=0,
ymax=4.5e-15,
ylabel style={font=\color{white!15!black}},
ylabel={\small $|E_k-E_0|$},
axis background/.style={fill=white},
xmajorgrids,
ymajorgrids,
legend style={at={(0.75,1.25)}, anchor=south west, legend cell align=left, align=left, draw=white!15!black,font=\small }
]
\addplot [color=red,dotted, line width=2.0pt]
  table[row sep=crcr]{%
0.00999999999999979	0\\
0.0200490524346391	8.88178419700125e-16\\
0.0301969791918193	8.88178419700125e-16\\
0.040495195628675	2.66453525910038e-15\\
0.0509975672957772	2.66453525910038e-15\\
0.0617613419063803	2.66453525910038e-15\\
0.0728481870032027	1.77635683940025e-15\\
0.0843253687500969	1.77635683940025e-15\\
0.0962671139563955	1.77635683940025e-15\\
0.108756206121519	1.77635683940025e-15\\
0.121885877269833	1.77635683940025e-15\\
0.135762070852591	8.88178419700125e-16\\
0.150506166965576	1.77635683940025e-15\\
0.166258278809961	1.77635683940025e-15\\
0.183181246415782	2.66453525910038e-15\\
0.201465464676261	2.66453525910038e-15\\
0.221334675904616	2.66453525910038e-15\\
0.243052808470375	2.66453525910038e-15\\
0.293340083457505	2.66453525910038e-15\\
0.322710610736245	3.5527136788005e-15\\
0.35554647503653	3.5527136788005e-15\\
0.392419828808678	3.5527136788005e-15\\
0.533555494766412	3.5527136788005e-15\\
0.592664363734248	3.5527136788005e-15\\
0.981929804073577	3.5527136788005e-15\\
1.07556664251516	3.5527136788005e-15\\
1.17328831059534	3.5527136788005e-15\\
1.27469553278405	3.5527136788005e-15\\
1.37945114779168	3.5527136788005e-15\\
1.48726565510495	3.5527136788005e-15\\
1.59788346565891	3.5527136788005e-15\\
1.71107180551196	3.5527136788005e-15\\
1.8266124652162	3.5527136788005e-15\\
2.18527436739457	3.5527136788005e-15\\
2.30816420961649	3.5527136788005e-15\\
2.43238399582543	3.5527136788005e-15\\
2.55772949236841	3.5527136788005e-15\\
2.68399504555366	3.5527136788005e-15\\
2.81097360264576	3.5527136788005e-15\\
3.06623572398324	3.5527136788005e-15\\
3.32184078813179	3.5527136788005e-15\\
3.44924740712205	3.5527136788005e-15\\
3.57611113682445	3.5527136788005e-15\\
4.19505022840563	3.5527136788005e-15\\
4.31433262774484	3.5527136788005e-15\\
4.43163975431987	3.5527136788005e-15\\
4.76962447236148	3.5527136788005e-15\\
4.87689967358183	3.5527136788005e-15\\
4.98106596282165	3.5527136788005e-15\\
5.08182621527071	3.5527136788005e-15\\
5.17883295250211	3.5527136788005e-15\\
5.27167422313305	3.5527136788005e-15\\
5.35986601975335	3.5527136788005e-15\\
5.44286368779437	3.5527136788005e-15\\
5.52010844186878	3.5527136788005e-15\\
5.59111677975761	3.5527136788005e-15\\
5.65559313374542	3.5527136788005e-15\\
5.71351558564297	3.5527136788005e-15\\
5.91986478491462	3.5527136788005e-15\\
5.94865827362874	4.44089209850063e-15\\
5.97457380049283	4.44089209850063e-15\\
5.9980320030886	3.5527136788005e-15\\
6.0389513448123	3.5527136788005e-15\\
6.05697105132563	3.5527136788005e-15\\
6.07366662212989	3.5527136788005e-15\\
6.08922301174809	3.5527136788005e-15\\
6.10379857103771	2.66453525910038e-15\\
6.11752971239285	8.88178419700125e-16\\
6.13053476593444	8.88178419700125e-16\\
6.14291716258854	2.66453525910038e-15\\
6.15476806717823	2.66453525910038e-15\\
6.16616856717082	2.66453525910038e-15\\
6.17719150524514	2.66453525910038e-15\\
6.18790302829445	1.77635683940025e-15\\
6.20863071311259	1.77635683940025e-15\\
};
\addlegendentry{$\text{Without VPA with tolerance=10}^{\text{-15}}$}

\addplot [color=blue, dashdotted, line width=2.0pt]
  table[row sep=crcr]{%
0.00999999999999979	0\\
0.0200490524346399	0\\
0.0301969791918211	0\\
0.0404951956286803	8.88178419700125e-16\\
0.0617613419063954	8.88178419700125e-16\\
0.0728481870032232	8.88178419700125e-16\\
0.35554647503675	8.88178419700125e-16\\
0.392419828808958	8.88178419700125e-16\\
0.730522839232465	8.88178419700125e-16\\
0.808859164734546	8.88178419700125e-16\\
1.48726565510775	8.88178419700125e-16\\
1.59788346566197	0\\
2.55772949237299	0\\
2.68399504555844	0\\
3.19410033150493	0\\
3.32184078813777	0\\
4.87689967358981	0\\
4.98106596282957	0\\
5.27167422314052	0\\
5.35986601976057	0\\
5.52010844187534	0\\
5.59111677976375	8.88178419700125e-16\\
5.97457380049555	8.88178419700125e-16\\
5.99803200309112	8.88178419700125e-16\\
6.19836391244434	8.88178419700125e-16\\
};
\addlegendentry{$\text{VPA with tolerance=10}^{\text{-15}}$}

\addplot [color=green, dashed, line width=2.0pt]
  table[row sep=crcr]{%
0.00999999999999979	0\\
0.0728481870032098	0\\
0.0843253687501058	0\\
0.293340083457468	0\\
0.322710610736207	0\\
5.99803200308934	0\\
6.01939033653293	0\\
6.19836391244326	0\\
};
\addlegendentry{$\text{VPA with tolerance=10}^{\text{-16}}$}

\addplot [color=black, dotted, line width=2.0pt]
  table[row sep=crcr]{%
0.00999999999999979	0\\
6.19836391244319	0\\
};
\addlegendentry{$\text{VPA with tolerance=10}^{\text{-17}}$}

\end{axis}

\begin{axis}[%
width=1.227\fheight,
height=0.658\fheight,
at={(-0.16\fheight,-0.072\fheight)},
scale only axis,
xmin=0,
xmax=1,
ymin=0,
ymax=1,
axis line style={draw=none},
ticks=none,
axis x line*=bottom,
axis y line*=left
]
\end{axis}
\end{tikzpicture}%
\caption{Energy error}
\label{fig: vpa07e}
    \end{subfigure}
    \hspace{1.5cm}
    \begin{subfigure}{.37\textwidth}
           \setlength\fheight{5.5 cm}
           \setlength\fwidth{\textwidth}
\raisebox{-62mm}{\input{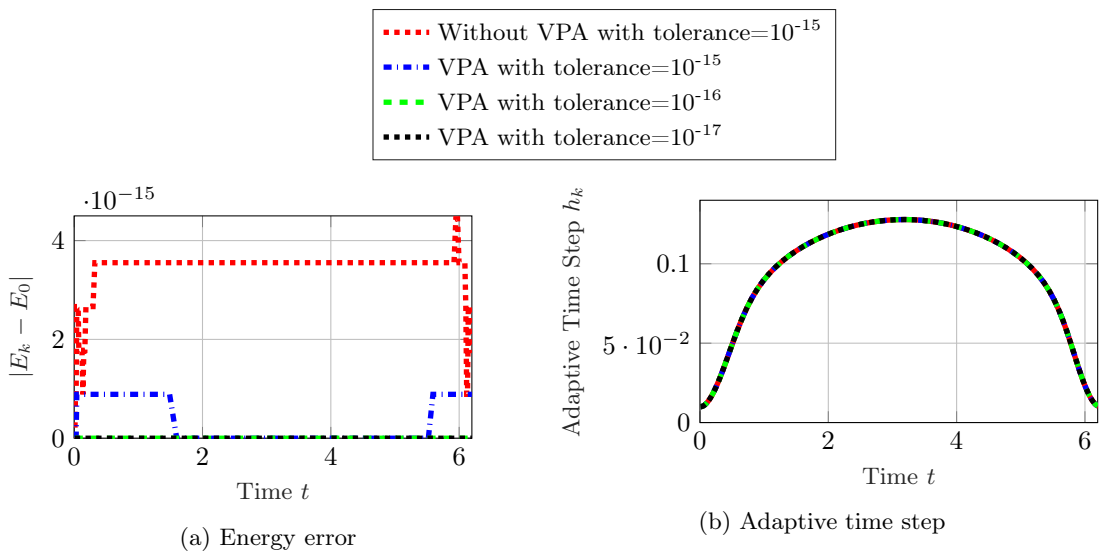}}
\caption{Adaptive time step}
\label{fig: vpa07h}
    \end{subfigure}    
    \caption{EpAVI algorithm for $e=0.7$ using VPA for different tolerance values.}
    \label{fig:vpa07}
\end{figure}
\par
It is important to note that the discrete energy error bound in \eqref{eq: dEbound} is very conservative because it assumes the worst case scenario where  the energy error is maximal and of the same sign in every step. We have studied Kepler's two-body problem using Variable Precision Arithmetic (VPA) to demonstrate that the EpAVI algorithm can resolve the issue of accumulation of discrete energy error by solving the discrete equations with higher precision. We have compared the numerical performance of the EpAVI algorithm using VPA with 18 significant digits for different tolerance values in Figure \ref{fig:vpa07}. The higher precision with VPA allows us to solve the discrete update equations with lower tolerance values. The energy error plots in Figure \ref{fig: vpa07e} show that the accumulation of discrete energy error decreases with decreasing tolerance, and VPA with $tolerance=10^{-17}$ achieves exact energy conservation up to precision. The adaptive time-step behavior in Figure \ref{fig: vpa07h} also matches our previous results from Section \ref{s:3}. 
\par 
Based on these results, we have shown that increasing the precision can lead to nearly exact energy behavior {\em without requiring time-step sizes to vanish}. It is important to note that using VPA makes the EpAVI algorithm computationally prohibitive. The key takeaway from these results obtained using VPA is that, if solved with enough precision, the adaptive algorithm can conserve the discrete energy exactly. In fact, the energy error results obtained by using VPA strongly indicate that the EpAVI algorithm is a {\em backward stable algorithm}. 
%
%

%
\section{Backward Error Analysis}
\label{s:6}
There are two challenges in showing backward stability for energy-preserving, adaptive time-step variational integrators. First, these algorithms are derived from the Lagrangian perspective and hence, the governing Euler-Lagrange equations are second-order.  The second challenge is the unique nature of the adaptive algorithm where the adaptive time-step is computed from the discrete Euler-Lagrange equations \eqref{eq:DEL}-\eqref{eq:DEE} and thus the time transformation map $t(a)$ in \eqref{eq:t(a)} is not explicitly given in the continuous-time setting. This is related to the fact that the Euler-Lagrange equation \eqref{eq:CEL} corresponding to the time variable $t$ is a redundant equation. 
\par
Vermeeren \cite{vermeeren2017modified} developed the theoretical framework for backward error analysis from the Lagrangian side and derived modified equations for fixed time-step variational integrators. Our goal is to apply backward stability results from \cite{vermeeren2017modified} to energy-preserving, adaptive time-step variational integrators. We first introduce a time transformation and show that the energy-preserving adaptive algorithm is equivalent to a fixed time-step variational integrator in the transformed setting. We then review some concepts required for the backward error analysis of variational integrators from a Lagrangian perspective and also modify some definitions to incorporate the time transformation. 
%
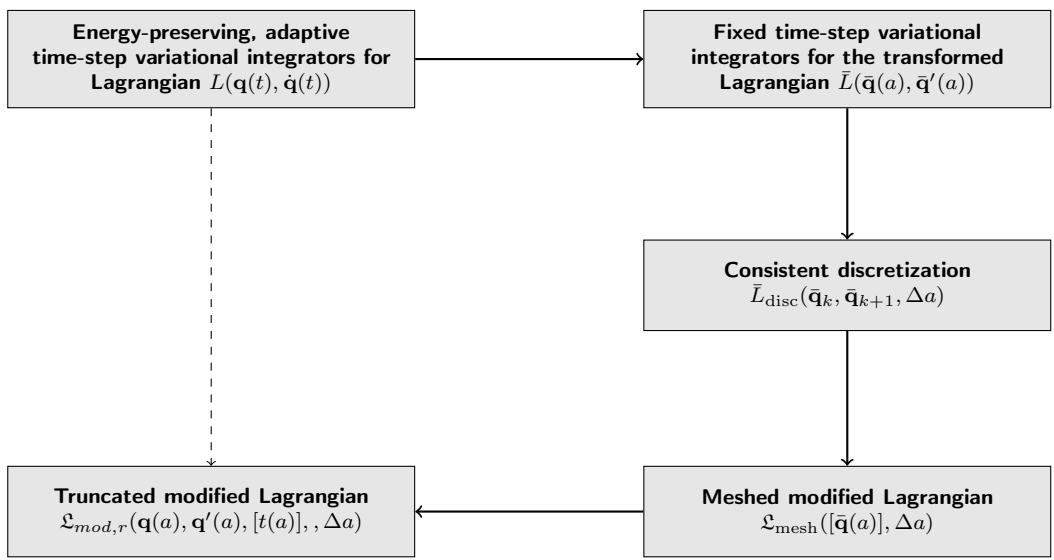
\begin{figure}[h]

\centering
\begin{tikzpicture}
[node distance = 1cm, auto,font=\footnotesize,
every node/.style={node distance=3cm},
comment/.style={rectangle, inner sep= 5pt, text width=5cm, node distance=0.25cm, font=\scriptsize\sffamily},
force/.style={rectangle, draw, fill=black!10, inner sep=5pt, text width=5cm, text badly centered, minimum height=1.2cm, font=\bfseries\footnotesize\sffamily}] 

\node [force] (VI) {Fixed time-step variational integrators for the transformed Lagrangian $\bar{L}(\qbar(a),\qbar'(a))$};
\node [force, left=3 cm of VI] (EpVI) {Energy-preserving, adaptive time-step variational integrators for Lagrangian $L(\q(t), \dot{\q}(t))$};
\node [force, below of=VI] (Emod) {Consistent discretization $\bar{L}_{\rm disc} (\qbar_k,\qbar_{k+1},\Delta a)$};
\node [force, below of=Emod] (Lmesh) {Meshed modified Lagrangian $\mathfrak L_{\rm mesh}([\qbar(a)],\Delta a)$};
\node [force, left=3 cm of Lmesh] (Lmod) {Truncated modified Lagrangian $\mathfrak L_{mod,r}(\q(a),\q'(a),[t(a)],,\Delta a)$};


;


;





\path[->,thick] 
(EpVI) edge (VI)
(VI) edge (Emod)
(Emod) edge (Lmesh)
(Lmesh) edge (Lmod);
\path[->,dashed] 
(EpVI) edge (Lmod);
\end{tikzpicture} 
\caption{Schematic overview of the backward error analysis of energy-preserving, adaptive time-step variational integrators.}
\label{fig:schematic}
\end{figure}
\subsection{Time Transformation}
\label{s:4}
In this section, we introduce a time transformation from the physical time $t$ to a fictitious time $a$ and represent the energy-preserving, adaptive time-step variational integrator (EpAVI) as a fixed time-step variational integrator in this transformed setting. 
Using ideas from extended Lagrangian mechanics, we interpret the parameterization of $t=t(a)$ in terms of the independent variable $a$ as a time transformation. For such a time transformation, i.e. $t \to a$, the configuration in the transformed setting is given by $\qbar(a)=(\q(a),t(a))$. For this new set of states, we introduce
\begin{equation}
     \bar{L}((\qbar(a),\qbar'(a)))=t'(a)L(\q(t),\qdot(t))=t'(a)L\left( \q(a),\frac{\q'(a)}{t'(a)}\right),
     \label{eq:t(a)}
\end{equation}
where we have used $\qdot(t)=\frac{\q'(a)}{t'(a)}$. The Euler-Lagrange equations in the transformed setting are 
\begin{align}
     \frac{\partial \bar{L}}{\partial \q} - \frac{d}{da} \left(  \frac{\partial \bar{L}}{\partial \q'} \right)&= t'(a)  \frac{\partial L}{\partial \q} - \frac{d}{da} \left( t'(a)  \frac{\partial L}{\partial \qdot} \frac{1}{t'(a)} \right)   = t'(a) \left( \frac{\partial L}{\partial \qdot} - \frac{d}{dt} \left(  \frac{\partial L}{\partial \qdot} \right)   \right)  =\bzero,
     \label{eq:CEL}\\
%
     \frac{\partial \bar{L}}{\partial t} - \frac{d}{da} \left(  \frac{\partial \bar{L}}{\partial t'} \right)  &= -\frac{d}{da} \left( L + t'(a) \frac{\partial L}{\partial \qdot} \left(-\frac{\qdot}{t'(a)}\right) \right) = t'(a) \frac{d}{dt} \left(  \frac{\partial L}{\partial \qdot} \cdot \dot{\q} - L \right) = 0.
        \label{eq:CEE}
\end{align}
Thus, the Euler-Lagrange equations in the transformed setting lead to the original Euler-Lagrange equation (\ref{eq:CEL}) along with the redundant energy equation (\ref{eq:CEE}). Although, these equations are the same as the extended Euler-Lagrange equations, it is important to interpret the above equations as Euler-Lagrange equations in the transformed setting. This interpretation will play a key role in the following Lemma about the equivalence of trajectories of the original Lagrangian system and the transformed system. \\
\\
\textbf{Lemma 1:} Let $\q(t)$ and $\q(a)$ be the trajectory solutions to the Euler-Lagrange equations corresponding to the physical Lagrangian $L$ and the transformed Lagrangian $\bar{L}$, respectively. Then, 
    \begin{equation*}
        \qbar(T)=(\q(\alpha(T)),\alpha(T)),
    \end{equation*}
where $\alpha(T)=\int_0^T t'(a) da$. 
\\
\textbf{Proof:} This follows from the fact that $\qbar(T)$ and $\q(\alpha(T))$ satisfy the original Euler-Lagrange equation with same initial conditions.
\par
Consider a discrete trajectory in the transformed setting with $\Delta a=a_{k+1}-a_k$. The discrete action $\bar{\mathfrak B}_d$ approximates the action integral
\begin{equation}
       \bar{\mathfrak B}_d =  \sum_{i=1}^N    \bar{L}_d(\qbar_k,\qbar_{k+1}) \approx \int_{a_0}^{a_f} \bar{L}(\qbar(a),\qbar'(a)) da, 
\end{equation}
where $\qbar_k=(\q_k,t_k)$ and the discrete Lagrangian is given by
\begin{equation}
         \bar{L}_d(\qbar_k,\qbar_{k+1})=\Delta a  \left[ \bar{L}\left(\frac{\qbar_k+ \qbar_{k+1}}{2},\frac{\qbar_{k+1}-\qbar_k}{\Delta a}\right)  \right] \approx \int_{a_k}^{a_{k+1}} \bar{L}(\qbar(a),\qbar'(a)) \ da .
\end{equation}
The resulting discrete Euler-Lagrange equations are given by
\begin{equation}
            \frac{\partial \bar{L}_d(\qbar_{k-1},\qbar_{k})}{\partial \q_k} + \frac{\partial \bar{L}_d(\qbar_k,\qbar_{k+1})}{\partial \q_k} =  \bzero,
            \label{eq:DEL}
\end{equation}
\begin{equation}
            \frac{\partial \bar{L}_d(\qbar_{k-1},\qbar_{k})}{\partial t_k} + \frac{\partial \bar{L}_d(\qbar_k,\qbar_{k+1})}{\partial t_k} =0.
            \label{eq:DEE}
\end{equation}
These variational integrators with a fixed time-step $\Delta a$ in the transformed setting are equivalent to the energy-preserving, adaptive time-step variational integrators. This equivalence will allow us to apply established results from backward error analysis of fixed time-step variational integrators.
\subsection{Consistent Discretization and Modified Equations}
\label{s:61}
Since variational integrators are written in terms of the discrete  Lagrangian $\bar L_d$, the first step is to define a consistent discretization of the continuous Lagrangian $\bar{L}(\qbar,\qbar')$. We define $\bar{L}_{\rm disc}$ as follows
\begin{equation}
   \bar{ \mathfrak{B}}_d=\sum_{k=0}^{N-1}\Delta a \bar L_{\rm disc}(\qbar_k,\qbar_{k+1},\Delta a),
\end{equation}
where $\bar{\mathfrak{B}}_d$ is the discrete action in the transformed setting. It is important to note that $\bar{L}_{disc}$ is different from the discrete Lagrangian $\bar  L_d$ defined earlier. In fact, $\bar L_d$ and $\bar L_{disc}$ are related by a nonzero scaling factor of $\Delta a$, i.e. $\bar L_d=\Delta a \bar L_{disc}$. A consistent discretization $\bar L_{disc}$ of $\bar L$ ensures that the discrete Euler-Lagrange equations \eqref{eq:DEL}-\eqref{eq:DEE}  are a consistent discretization of continuous Euler-Lagrange equations \eqref{eq:CEL}-\eqref{eq:CEE}.  For example, the midpoint rule used in this work is an example of consistent discretization with 
    \begin{equation*}
         \bar{L}_{\rm disc}= \left[ \bar{L}\left(\frac{\qbar_k+ \qbar_{k+1}}{2},\frac{\qbar_{k+1}-\qbar_k}{\Delta a}\right)  \right] .
        \end{equation*}
As an example, consider the Kepler's two-body problem from Section \ref{s:3}. The governing discrete equations based on the midpoint rule
\begin{equation}
\left(\frac{q^1_{k+1}-q^1_k}{t_{k+1}-t_k} \right) - (t_{k+1}-t_k) \left( \frac{q^1_k + q^1_{k+1}}{4\left(\left(q^1_k + q^1_{k+1}\right)^2 +\left(q^2_k + q^2_{k+1}\right)^2\right)^{\frac{3}{2}}}\right)=p^1_k,
\end{equation}
\begin{equation}
\left(\frac{q^2_{k+1}-q^2_k}{t_{k+1}-t_k} \right) - (t_{k+1}-t_k) \left( \frac{q^2_k + q^2_{k+1}}{4\left(\left(q^1_k + q^1_{k+1}\right)^2 +\left(q^2_k + q^2_{k+1}\right)^2\right)^{\frac{3}{2}}}\right)=p^2_k,
\end{equation}
\begin{equation}
\frac{1}{2} \left( \frac{q^1_{k+1}-q^1_k}{t_{k+1}-t_k} \right)^2 + \frac{1}{2} \left( \frac{q^2_{k+1}-q^2_k}{t_{k+1}-t_k} \right)^2 - \frac{2}{\sqrt{\left(q^1_k + q^1_{k+1}\right)^2 +\left(q^2_k + q^2_{k+1}\right)^2}}=E_k.
\end{equation}
These are a consistent discretization of the continuous Euler-Lagrange equations.
\par
Since the governing equations from the Lagrangian perspective, i.e. Euler-Lagrange equations, are second-order equations, the modified equation for the discrete Euler-Lagrange equations should also be second-order. As discussed in Section \ref{s:2}, the energy-preserving, adaptive time-step variational integrators are based on the extended Lagrangian mechanics and its discrete analogue.  To use tools from standard backward analysis of numerical integrators, we first need to define what constitutes a modified equation in the transformed setting.  We define $ \pmb{\psi}=[\pmb{\psi}_{EL},\psi_E]$  where $\pmb \psi_{EL}$ is the discrete Euler-Lagrange equation
    \begin{equation}
        \pmb \psi_{EL}(\qbar_{k-1},\qbar_{k},\qbar_{k+1},\Delta a) =  \frac{\partial \bar{L}_d(\qbar_{k-1},\qbar_{k})}{\partial {\q}_k} + \frac{\partial \bar{L}_d(\qbar_k,\qbar_{k+1})}{\partial {\q}_k} =  \bzero,
    \end{equation}
and $\psi_E$ is the discrete energy equation 
     \begin{equation}
        \psi_E(\qbar_{k-1},\qbar_{k},\qbar_{k+1},\Delta a) =  \frac{\partial \bar{L}_d(\qbar_{k-1},\qbar_{k})}{\partial t_k} + \frac{\partial \bar{L}_d(\qbar_k,\qbar_{k+1})}{\partial t_k} =  0. 
    \end{equation}
\textbf{Definition:} The differential equation $\qbar''=\mathbf f(\qbar(a),\qbar'(a),\Delta a)$ is a modified equation for the discrete equation  $\psi(\qbar_{k-1},\qbar_{k},\qbar_{k+1},\Delta a)=\bzero$ if, for every $r$, every admissible family $(\qbar_{\Delta a})$ of solutions of 
\begin{equation} 
\qbar''_{\Delta a}=\tau_r  (\mathbf f(\qbar_{\Delta a},\qbar'_{\Delta a},\Delta a)) ,
\end{equation}
where $\tau_r$ denotes truncation after the $\Delta a^r$- term, satisfies
\begin{equation*}
   \pmb  \psi(\qbar_{\Delta a}(a-\Delta a),\qbar_{\Delta a}(a),\qbar_{\Delta a}(a+\Delta a),\Delta a ) = O(\Delta a^{r+1}),
\end{equation*}
as $\Delta a\to 0$ for all $a$. 
 \begin{remark}
 A family of curves $(\qbar_{\Delta a})_{\Delta a \in (0,\infty)}$ of smooth curves $\qbar_{\Delta a}:[a_h,b_h] \to \Real^N$ is called admissible if there exists $\Delta a_{max}>0$ such that for each $k\geq 0,||\qbar_{\Delta a}^{(k)}||_{\infty} $ is bounded as a function of $\Delta a \in (0,\Delta a_{max}]$, where $|| \cdot ||_{\infty}$ denotes the supremum norm. These admissible families are families of curves parametrized by $\Delta a$ whose derivatives do not blow up when $\Delta a \to 0$. Admissibility of a family of smooth curves guarantees that in power series expansions like $\qbar_{\Delta a}(a+\Delta a)=\qbar_{\Delta a}(a) + \Delta a \qbar'_{\Delta a}(a) + \cdots$ the asymptotic behavior of each term is determined by the exponent of $\Delta a$ in that term.
 \end{remark}
\subsection{Modified Lagrangian}
\label{s:62}
For a variational integrator in the transformed setting with a modified equation $\pmb \psi$, ideally we would like to construct a modified Lagrangian $\bar L_{\rm mod}$ such that modified equation $\pmb \psi$ is found through the Euler-Lagrange equations. However, the modified equations generally have the form of non-converging power series so we will focus on the $r-$truncation of the modified equation. Thus for the backward error analysis, our goal is to obtain a modified Lagrangian $\bar{L}_{\rm mod}$ such that the $r-$truncation of the Euler-Lagrange equations for the modified Lagrangian are the $r-$truncation of the modified equations corresponding to the energy-preserving, adaptive algorithm. 
\par
If $\qbar$ is a curve such that $(\qbar(a_0 + k\Delta a))_{k\in [0,\cdots, n]}$ is critical for the discrete action $\bar{\mathfrak{B}}_d$, then $\qbar$ solves the meshed variational problem for $\bar{\mathfrak L}_{\rm mesh}$
    \begin{equation}
      \bar{  \mathfrak{B}}_d((\qbar(a_0 + k\Delta a)))= \sum_{k}\bar{L}_d(\qbar_k,\qbar_{k+1},\Delta a) \simeq \int^{a_0 + n\Delta a}_{a_0}\bar{\mathfrak L}_{\rm mesh} ([\qbar(a)],\Delta a) \ da , 
    \end{equation}
    where $[\qbar(a)]$ contains $\qbar(a)$ and its derivatives w.r.t fictitious time $a$. The meshed modified Lagrangian $\bar{\mathfrak L}_{\rm mesh}$ is defined by the following formal power series
\begin{equation}
    \bar{\mathfrak L}_{\rm mesh}([\qbar],\Delta a):=  \bar{\mathfrak L}_{\rm disc}([\qbar],\Delta a) -\frac{\Delta a^2}{24}\frac{d^2}{dt^2}\bar{\mathfrak L}_{\rm disc}([\qbar],\Delta a) + \frac{7\Delta a^4}{5760}\frac{d^4}{dt^4}\bar{\mathfrak L}_{\rm disc}([\qbar],\Delta a) + \cdots  \ .
    \label{eq:Lmeshmod}
\end{equation}
Since $\bar{\mathfrak L}_{\rm disc}([\qbar],\Delta a)$ is a consistent discretization of $\bar L$, it is clear that zeroth order term of the meshed modified is the original continuous Lagrangian, i.e. $\bar {\mathfrak L}_{\rm mesh}([\qbar],\Delta a)=\bar L(\qbar, \qbar') + O(\Delta a)$.
\begin{remark}
It is evident from \eqref{eq:Lmeshmod} that the meshed modified Lagrangian $\bar{\mathfrak L}_{mesh}$ depends on higher derivatives of $\qbar$ instead of just $\qbar$ and $\qbar'$. This is related to the unconventional nature of the meshed variational problem that defines $\bar{\mathfrak L}_{mesh}$. Unlike the classical variational problem that seeks critical curves of the action integral \eqref{eq:action} within the set of smooth curves $\mathcal C^{\infty}$, the meshed variational problem with mesh size $\Delta a$ seeks smooth curves that for every $a_0 \in \Real$ are critical for the action integral within the set of piecewise smooth curves $\mathcal M^{a_0,\Delta a}$ whose nonsmooth points lie in the mesh $a_0 + \Delta a \mathbb{Z}$. For a thorough exposition, the interested reader may consult
 Section 4 in \cite{vermeeren2017modified}.
\end{remark}
\par 
Based on Lemma $8$ from \cite{vermeeren2017modified}, the modified equation for variational integrators can be calculated from the Euler-Lagrange equations of the meshed modified Lagrangian. Due to the presence of higher order derivatives in $\bar {\mathfrak L}_{\rm mesh}$, these modified equations are written as a power series in the fixed time-step $\Delta a$ and are usually non-convergent. In such cases, we are often only interested in truncated power series expressions so we also need to relax our notion of critical curves for both continuous and discrete variational problems. Consider the following definitions: 
\begin{enumerate}
\item \textbf{Classical $r$--critical}: An admissible family $(\qbar_{\Delta a})_{\Delta a \in (0, \infty)}$ of curves $\qbar_{\Delta a}:[0,1] \to \Real$ is $r$--critical for a family of actions $\bar{\mathfrak{B}}$ if every family of smooth variations $\delta \qbar_{\Delta a} $ satisfies $\delta \bar{\mathfrak{B}}_{\Delta a}=O(\Delta a^{r+1})|| \delta \qbar_{\Delta a} ||$. An admissible family $(\qbar_{\Delta a})_{\Delta a \in(0,\infty)}$ of curves $\qbar_{\Delta a}:[0,1] \to \Real$ is $r-$critical for a family of actions $\bar{\mathfrak{B}}_{\Delta a}=\int_0^1\bar{L}_{\Delta a} dt$ if and only if
\begin{equation}
    \bigg| \bigg| \frac{\delta \bar{L}_{\Delta a}}{\delta \qbar_{\Delta a}} \bigg|\bigg|_{\infty}=O(\Delta a^{r+1}) . 
\end{equation}
\item \textbf{Meshed $r$--critical}: An admissible family $(\qbar_{\Delta a})_{\Delta a \in (0, \infty)}$ of curves $\qbar_{\Delta a}:[0,1] \to \mathbb{R}$ is  meshed $r$--critical for a family of actions $\bar{\mathfrak{B}}$ if for every family of smooth variations $\delta \qbar_{\Delta a} \in \mathcal{M}^{a_{\Delta a},\Delta a} $ (i.e. piecewise smooth $\delta \qbar$ with nonsmooth points in a mesh of size $\Delta a$) we have $\delta \bar{\mathfrak{B}}_{\Delta a}=O(\Delta a^{r+1})|| \delta \qbar_{\Delta a} ||$. Similarly, an admissible family  $(\qbar_{\Delta a})_{\Delta a \in(0,\infty)}$ of curves $\qbar_{\Delta a}:[0,1] \to \mathbb{R}$ is meshed  $r$--critical for a family of actions $\bar{\mathfrak{B}}_{\Delta a}=\int_0^1\bar{L}_{\Delta a} dt$  if and only if
  \begin{equation}
       \bigg| \bigg|  \frac{\delta \bar{L}_{\Delta a}}{\delta \qbar_{\Delta a}} \bigg|\bigg|_{\infty}=O(\Delta a^{r+1}) \quad \mbox{and} \quad  \bigg| \bigg|  \frac{\delta \bar{L}_{\Delta a}}{\delta \qbar_{\Delta a}^{(l)}} \bigg|\bigg|_{\infty}=O(\Delta a^{r+l+1}) \ \mbox{for each} \ l \geq 2 . 
  \end{equation}
\item \textbf{Discrete $r$--critical}: A family $(\qbar_{\Delta a})_{\Delta a \in (0, \infty)}$ of discrete curves $\qbar_{\Delta a}:[0,1] \to \mathbb{R}$ is $r-$critical for a family of discrete actions $\bar{\mathfrak{B}}_d$ if for every family of discrete variations $\delta \qbar_{\Delta a} $ we have $\delta \bar{\mathfrak{B}}_{d, \Delta a}=O(\Delta a^{r+1})|| \delta \qbar_{\Delta a} ||$. A family of discrete curves $\qbar_{\Delta a}=(\qbar_{\Delta a,k})_{k \in [0, \cdots,n_k]}$ is $r$--critical for a family of discrete actions $\bar{\mathfrak{B}}_{d,\Delta a}=\sum \bar{L}_d$ if and only if 
  \begin{equation}
      \frac{\partial  \bar{L}_d}{\partial \qbar_{\Delta a,k}}(\qbar_{\Delta a,k},\qbar_{\Delta a,k+1}) + \frac{\partial  \bar{L}_d}{\partial \qbar_{\Delta a,k}}(\qbar_{\Delta a,k-1},\qbar_{\Delta a,k})=O(\Delta a^{r+1}) . 
  \end{equation}

\end{enumerate}
\par 
The $r-$critical families of curves for the meshed modified Lagrangian $\bar{\mathfrak{L}}_{\rm mesh}$ satisfy
\begin{equation}
\frac{\partial \bar{\mathfrak{L}}_{mesh}}{\partial {\q}}  - \frac{d}{da} \left( \frac{\partial \bar{\mathfrak L}_{\rm mesh}}{\partial \q'} \right) = O(\Delta a^{r+1}).
\end{equation}
Based on Proposition $6$ from \cite{vermeeren2017modified}, the modified equation can be written in the following form
\begin{equation}
\epsilon_0(\qbar, \qbar', \qbar'') + \Delta a \epsilon_1(\qbar, \qbar', \qbar'') + \Delta a^2{\epsilon}_2(\qbar, \qbar', \qbar'',\qbar^{(3)}) + \cdots + \Delta a^r {\epsilon}_r(\qbar, \qbar', \cdots ,\qbar^{(r+1)})= O(\Delta a^{r+1}).
\label{eq:recursion}
\end{equation}
where for $r \geq 1$, $\epsilon_r$ depends on $\qbar, \qbar', \cdots, \qbar^{(r+1)}$ but not on higher derivatives of $\qbar$. Using $\bar {\mathfrak L}_{\rm mesh}([\qbar],\Delta a)=\bar L(\qbar, \qbar') + O(\Delta a)$, it is easy to show that $\epsilon_0(\qbar, \qbar',\qbar'')$ is simply the continuous Euler-Lagrange equation \eqref{eq:CEL}. For a consistent discretization $\bar L_{\rm disc}$ and sufficiently small $\Delta a$, we can solve $\epsilon_0(\qbar, \qbar', \qbar'') + \Delta a \epsilon_1(\qbar, \qbar', \qbar'')=O(\Delta a^2)$ for $\q''$, i.e.
 \begin{equation}
\q''=\pmb F_{1}(\q,\q',[t],\Delta a) + O(\Delta a^{2}). 
\end{equation}
Following the recursion proposed in Section 4.6 of \cite{vermeeren2017modified}, we can show that the modified equation \eqref{eq:recursion} can be solved for $\q''$ in terms of $(\q,\q',[t],\Delta a)$
\begin{equation}
\q''=\pmb F_{r}(\q,\q',[t],\Delta a) + O(\Delta a^{r+1}). 
\end{equation}
It is important to note that we can only solve the modified equation for $\q$ and not for $\qbar$.  In other words, we can solve \eqref{eq:recursion} for $\q''$ in terms of $\q,\q'$ and arbitrary number of derivatives of $t$.  However, the redundant nature of \eqref{eq:CEE}  prevents us from solving the modified equations for $\qbar=[\q,t]$ jointly for $\qbar''$ in terms of $\qbar$ and $\qbar'$. 
\begin{remark} 
For the fixed time-step variational integrator considered in  \cite{vermeeren2017modified}, the author used the properties of admissible curves to write the modified equation as a second order differential equation with $O(\Delta a^{r+1})$ defect. This was achieved by a recursive process in Section 4.6 of \cite{vermeeren2017modified} where the leading order terms in the modified equation are solved for $\ddot{\q}$. This result doesn't carry over to our transformed setting due to the redundant nature of the continuous Euler-Lagrange equation \eqref{eq:CEE} corresponding to the time variable.  
\end{remark}
\par
Unlike the classical Lagrangian defined on the tangent bundle, the meshed modified Lagrangian $\bar{\mathfrak L}_{\rm mesh}$ also depends on higher derivatives of $\qbar(a)$. From the meshed modified Lagrangian $\bar {\mathfrak L}_{\rm mesh}$, we define the classical modified Lagrangian $\bar{\mathfrak L}_{\rm mod}$
\begin{equation}
\bar{\mathfrak L}_{\rm mod}(\q,\q',[t],\Delta a)=\bar{\mathfrak L}_{\rm mesh}([\qbar],\Delta a)\bigg|_{\q''=\pmb \psi(\q, \q',[t],\Delta a), \q'''= \frac{d}{da} \pmb \psi( \q, \q',[t],\Delta a), \cdots  },
\end{equation}
where $\q''=\pmb \psi(\q,\q',[t],\Delta a)$ is the modified equation. The $r-$th truncation of the modified Lagrangian is denoted by $\bar {\mathfrak L}_{\rm mod,r}$,
\begin{equation}
\bar{\mathfrak L}_{\rm mod,r}({\q}, \q',[t],\Delta a)=\tau_r\left( \bar{\mathfrak L}_{\rm mod}( \q,\q',[t],\Delta a) \right)=\tau_r\left(\bar{\mathfrak L}_{\rm mesh}([\qbar],\Delta a)\bigg|_{\q^{(r)}=F^{(j)}_{r-2}( \q, \q',[t],\Delta a)}  \right).
\end{equation} 
From the definition of $r-$th truncation of  modified Lagrangian, it follows that $\bar {\mathfrak L}_{\rm mod,r}({\q},{\q}',[t], \Delta a)=\bar {\mathfrak L}_{\rm mesh}( [\qbar] , \Delta a) + O(\Delta a^{r+1})$ for meshed $r-$critical families of curves. In fact, Vermeeren \cite{vermeeren2017modified} showed that the meshed modified Lagrangian $\bar {\mathfrak L}_{\rm mesh}$ and truncated modified Lagrangian $\bar{\mathfrak L}_{\rm mod,r}$ have the same $r-$critical families of curves. 
\begin{theorem} For a discretized Lagrangian $\bar{L}_{\rm disc}$ that is a consistent discretization of a regular Lagrangian $\bar L(\qbar(a),\qbar'(a))$, the $r-$th truncation of the Euler-Lagrange equation of $\bar{\mathfrak L}_{\rm mod,r}({\q}(a),{\q}'(a),[t(a)],\Delta a)$, solved for ${\q}''(a)$, is the $r-$th truncation of the modified equation. 
\end{theorem}
\begin{proof} This result follows directly from applying backward stability theorem from \cite{vermeeren2017modified} to the transformed Lagrangian system defined in Section \ref{s:4}. 
\end{proof}
Although our theorem is based on the application of ideas from \cite{vermeeren2017modified} to the transformed setting, the final result for the energy-preserving algorithm has two key differences due to the redundant nature of the energy equation.
\begin{enumerate}
\item The state vector $\qbar$ in the transformed setting contains both $\qbar$ and $t$, the modified equation is only obtained and solved for $\q''$. This is due to the fact that the modified equation corresponding to the $t$ variable is redundant.
\item The modified Lagrangian $\bar{\mathfrak L}_{\rm mod,r}$ in the transformed setting depends on $\q,\q',\Delta a$ and arbitrary high derivatives of $t$. This again is related to the redundant nature of the energy equation which prevents us from solving the modified equations jointly for $\qbar''$.
\end{enumerate}
\subsection{Mechanical System Example}
\label{s:63}
We consider a one degree of freedom mechanical system with a separable Lagrangian to illustrate the backward stability result shown in Section \ref{s:62}. Consider a Lagrangian system $L(q,\dot{q})=\frac{1}{2}m\dot{q}^2 - V(q)$ where $V(q)$ is the potential energy of the mechanical system. We employ the time transformation \eqref{eq:t(a)} to obtain the transformed Lagrangian
\begin{equation}
\bar L(\bar q(a), \bar q'(a))=t'(a)L\left(q(a), \frac{q'(a)}{t'(a)}\right)=t'(a) \left( \frac{1}{2}m\left( \frac{q'(a)}{t'(a)} \right)^2 - V(q(a)) \right).
\end{equation}
We use the midpoint rule to obtain the discretized Lagrangian $\bar L_{disc}$
\begin{equation}
\bar L_{\rm disc}(\bar q_k, \bar q_{k+1},\Delta a)=\frac{t_{k+1}-t_k}{\Delta a} \left( \frac{1}{2}m\left( \frac{q_{k+1}-q_k}{t_{k+1}-t_k} \right)^2 - V\left(\frac{q_{k+1}+q_k}{2}\right) \right).
\end{equation}
Using the power series expansions for $q_{k\pm1}(a)=q(a) \pm \frac{\Delta a}{2}q'(a) + \frac{\Delta a^2}{8}q''(a) \pm \frac{\Delta a^3}{24}q'''(a) + O(\Delta a^4) $ and $t_{k\pm1}(a)=t(a) \pm \frac{\Delta a}{2}t'(a) + \frac{\Delta a^2}{8}t''(a) \pm \frac{\Delta a^3}{24}t'''(a) + O(\Delta a^4) $, we obtain
\begin{align}\nonumber
\bar{\mathfrak L}_{\rm disc}([\bar q],\Delta a)&=\left( t' + \frac{\Delta a^2}{24}t''' + \cdots \right) \left( \frac{1}{2}m\left( \frac{ q' + \frac{\Delta a^2}{24}q''' + \cdots }{t' + \frac{\Delta a^2}{24}t''' + \cdots } \right)^2  - V\left(q + \frac{\Delta a^2}{8}q'' + \cdots \right) \right) \\ 
&=t'\left( \frac{1}{2}m \left( \frac{q'}{t'} \right)^2 -V \right)  + \frac{\Delta a^2}{24} \left( \frac{mq'q'''}{t'}- \frac{m(q')^2t'''}{2(t')^2} - 3q'' t'V_q - t'''V \right) + O(\Delta a^4),
\end{align}
where $V$ with no arguments is $V(q(a))$ and $V_q=\frac{\partial V}{\partial q}(q(a))$. From $\bar{\mathfrak{L}}_{\rm disc}([\bar q],\Delta a)$, we calculate the meshed modified Lagrangian in the transformed setting 
\begin{align}\nonumber
\bar{\mathfrak L}_{\rm mesh}([\bar q],\Delta a)&= \bar{\mathfrak L}_{\rm disc}([\bar q],\Delta a) - \frac{\Delta a^2}{24} \frac{d^2}{da^2} \bar{\mathfrak L}_{\rm disc}([\bar q],\Delta a)  + O(\Delta a^4) \\
&=t'\left( \frac{1}{2}m \left( \frac{q'}{t'} \right)^2 -V \right)  + \frac{\Delta a^2}{24} \left(-\frac{m(q'')^2}{t'}  + \frac{2mq'q''t''}{(t')^2} - \frac{m(q')^2(t'')^2}{(t')^3} \right) \nonumber \\ 
&+ \frac{\Delta a^2}{24} \left( 2q't''V_q + (q')^2t'V_{qq} - 2q'' t'V_q \right) + O(\Delta a^4). 
\label{eq:Lmesh}
\end{align}
The modified Euler-Lagrange and modified energy equation up to second order in the transformed setting are obtained from
\begin{align}\nonumber
O(\Delta a^4)&= \frac{\partial \bar {\mathfrak L}_{\rm mesh}}{\partial q} - \frac{d}{da}\left(\frac{\partial \bar{\mathfrak L}_{\rm mesh}}{\partial q'}\right) \\
&= -\frac{mq''}{t'}+\frac{mq't''}{(t')^2}  -t'V_q + \frac{\Delta a^2}{24} \left(- (q')^2t'V_{qqq}-2t'''V_q -2q't''V_{qq}   - 4q''t'V_{qq} \right) \nonumber \\ 
&-   \frac{\Delta a^2}{24} \left(  \frac{2mq'''t''}{(t')^2} +  \frac{2mq''t'''}{(t')^2} -  \frac{6mq''(t'')^2}{(t')^3} - \frac{4mq't''t'''}{(t')^3}+ \frac{6mq'(t'')^3}{(t')^4} \right), \\\nonumber
O(\Delta a^4)&= \frac{\partial \bar {\mathfrak L}_{mesh}}{\partial t} - \frac{d}{da}\left(\frac{\partial \bar{\mathfrak L}_{mesh}}{\partial t'}\right)\\
&=\left(\frac{q'}{t'} \right)\left(-\frac{mq''}{t'} + \frac{mq't''}{(t')^2} - t'V_q \right) + O(\Delta a^2),
\end{align}
where $V_{qq}=\frac{\partial^2 V}{\partial q^2}(q(a))$ and $V_{qqq}=\frac{\partial^3 V}{\partial q^3}(q(a))$. 
\par 
Due to the redundant nature of the leading order in modified energy equation, the modified equations for the transformed system can not be solved recursively for $\bar q''=[q''(a), t''(a)]$. This is a consequence of the redundancy in the Euler-Lagrange equations in the transformed setting. We solve recursively for $q''$ in terms of $(q,q')$ and $(t',t'',t''')$. In the leading order we have $q''=\left( \frac{q't''}{t'} - \frac{t'^2V_q}{m} \right)$. Substituting this in the second-order terms in the modified Euler-Lagrange equation, we get
\begin{align} \nonumber
0
&= \frac{\partial \bar {\mathfrak L}_{\rm mesh}}{\partial q} - \frac{d}{da}\left(\frac{\partial \bar{\mathfrak L}_{\rm mesh}}{\partial q'}\right) + O(\Delta a^4) \\ \nonumber
&= -\frac{mq''}{t'}+\frac{mq't''}{(t')^2}  -t'V_q + \frac{\Delta a^2}{24} \left(- (q')^2t'V_{qqq}-2t'''V_q -2q't''V_{qq}   - 4\left( \frac{q't''}{t'} - \frac{(t')^2V_q}{m} \right)t'V_{qq} \right) \\ \nonumber
&-   \frac{\Delta a^2}{24} \left(  \frac{2m\frac{d}{da}\left( \frac{q't''}{t'} - \frac{(t')^2V_q}{m} \right)t''}{(t')^2} +  \frac{2m\left( \frac{q't''}{t'} - \frac{(t')^2V_q}{m} \right)t'''}{(t')^2} -  \frac{6m\left( \frac{q't''}{t'} - \frac{t'^2V_q}{m} \right)(t'')^2}{(t')^3} \right) \\ \nonumber
&+   \frac{\Delta a^2}{24} \left(  \frac{4mq't''t'''}{(t')^3}- \frac{6mq'(t'')^3}{(t')^4} \right) + O(\Delta a^4)  \\
&= -\frac{mq''}{t'}+\frac{mq't''}{(t')^2}  -t'V_q + \frac{\Delta a^2}{24} \left( \frac{4(t')^3V_qV_{qq}}{m} -4q't''V_{qq}- (q')^2t'V_{qqq}  \right) + O(\Delta a^4), 
\end{align}
which can be re-written as
\begin{equation}
q''= \frac{q't''}{t'}  -\frac{(t')^2V_q}{m} + \frac{\Delta a^2}{24m} \left( \frac{4(t')^4V_qV_{qq}}{m} -4q't't''V_{qq}- (q')^2(t')^2V_{qqq}\right) + O(\Delta a^4).
\label{eq:EL_mesh}
\end{equation}
We also obtain the modified Lagrangian up to second order 
\begin{align} \nonumber
\bar{\mathfrak L}_{\rm mod}& =\bar{\mathfrak L}_{\rm mesh}([\bar q],\Delta a)\bigg|_{q''=\left( \frac{q't''}{t'} - \frac{t'^2V_q}{m} \right)} \\ \nonumber
&=t'\left( \frac{1}{2}m \left( \frac{q'}{t'} \right)^2 -V \right)  + \frac{\Delta a^2}{24} \left( \frac{2mq'\left( \frac{q't''}{t'} - \frac{t'^2V_q}{m} \right)t''}{(t')^2}-\frac{m\left( \frac{q't''}{t'} - \frac{(t')^2V_q}{m} \right)^2}{t'}   - \frac{m(q')^2(t'')^2}{(t')^3} \right) \\ \nonumber
&+ \frac{\Delta a^2}{24} \left( 2q't''V_q + (q')^2t'V_{qq} - 2\left( \frac{q't''}{t'} - \frac{(t')^2V_q}{m} \right) t'V_q \right) + O(\Delta a^4)\\
&=t'\left( \frac{1}{2}m \left( \frac{q'}{t'} \right)^2 -V \right)  + \frac{\Delta a^2}{24} \left( \frac{t'^3(V_q)^2}{m} + q'^2t'V_{qq} \right) +  O(\Delta a^4).
\end{align}
Unlike the meshed modified Lagrangian \eqref{eq:Lmesh} which depends on higher derivatives of $q$ and $t$, the modified Lagrangian in this case only depends on $\bar{q}=(q,t)$ and first derivative $\bar{q}'=(q',t')$.
The Euler-Lagrange equation of $\bar{\mathfrak L}_{\rm mod,3}$ is
\begin{align}\nonumber
0&= \frac{\partial \bar {\mathfrak L}_{\rm mod,3}}{\partial q} - \frac{d}{da}\left(\frac{\partial \bar{\mathfrak L}_{\rm mod,3}}{\partial q'}\right)\\\nonumber
&=-t'V_q + \frac{\Delta a^2}{24} \left( \frac{2t'^3V_qV_{qq}}{m} + V_{qqq}q'^2t' \right) - \frac{d}{da}\left( \frac{mq'}{t'} + \frac{\Delta a^2}{24} \left( 2V_{qq}q't' \right) \right) \\
&=-\frac{mq''}{t'}+\frac{mq't''}{(t')^2} -t'V_q + \frac{\Delta a^2}{24} \left( \frac{2t'^3V_qV_{qq}}{m}   - 2q'V_{qq} t'' - V_{qqq}q'^2t' -2V''q''t'  \right).
\label{eq:EL_mod}
\end{align}
It is important to note that the Euler-Lagrange equation for $\bar{\mathfrak L}_{\rm mod,3}$ does not contain an error term. However, when we solve the Euler-Lagrange equation \eqref{eq:EL_mod} for $q''$
\begin{align}\nonumber
q''&= \frac{\frac{mq't''}{(t')^2} -t'V_q + \frac{\Delta a^2}{24} \left( \frac{2(t')^3V_qV_{qq}}{m}   - 2q't''V_{qq} - q'^2t'V_{qqq} \right)}{\left( \frac{m}{t'} + \frac{\Delta a^2 t'V_{qq}}{12} \right)}\\\nonumber
&=\left[ \frac{q't''}{t'} -\frac{(t')^2V_q}{m} + \frac{\Delta a^2}{24m} \left( \frac{2(t')^4V_qV_{qq}}{m}   - 2q't't''V_{qq} - (q't')^2V_{qqq} \right) \right]\left( 1 - \Delta a^2  \frac{(t')^2V_{qq}}{12m} + O(\Delta a^4) \right)\\
&=\frac{q't''}{t'} -\frac{(t')^2V_q}{m} + \frac{\Delta a^2}{24m} \left( \frac{4(t')^4V_qV_{qq}}{m}   - 4q't't''V_{qq} - (q')^2(t')^2V_{qqq} \right) + O(\Delta a^4), 
\end{align}
we again get the modified equation \eqref{eq:EL_mesh}. Thus, the modified equation for the variational integrator based on the midpoint rule, solved for $q''$, is $O(\Delta a^4)$ close to the Euler-Lagrange equation of the truncated modified Lagrangian $\bar{\mathfrak L}_{\rm mod,3}$.
\section{Conclusions}
\label{s:7}
We have presented a numerical and theoretical assessment of energy-preserving, adaptive time-step variational integrators. We first compared the time adaptation and energy performance of the energy-preserving adaptive algorithm with the adaptive variational integrators for Kepler's two-body problem. We then obtained a conservative bound on the discrete energy error using the discrete energy equation and implemented the energy-preserving adaptive algorithm using VPA to demonstrate exact discrete energy conservation up to machine precision. 
\par 
Motivated by these numerical results, we proposed a time transformation approach to investigate the numerical stability of the adaptive algorithm. Using the equivalence of trajectories, we derived the energy-preserving, adaptive time-step variational integrators as fixed time-step variational integrators in the transformed setting. We then applied results from \cite{vermeeren2017modified} in the transformed setting to show that the $r-$th truncation of the modified equation for the variational integrator in the transformed setting, solved for $\q''$, is $O(\Delta a^{r+1})$ close to the Euler-Lagrange equation of the $r-$truncated modified Lagrangian $\bar{\mathfrak L}_{\rm mod,r}$. Finally, we constructed a modified Lagrangian for a midpoint rule discretization of a simple mechanical system.
\par
For future work, we would like to extend the backward stability results from this paper to the energy-preserving, adaptive time-step Lie group variational integrators. We would also like to study the connections between modified equations of adaptive variational integrators based on a prescribed monitor function and energy-preserving, adaptive time-step variational integrators.
\\ \\
\noindent
\textbf{Declaration of interest:} None.
\\ 
\\
\textbf{Funding:} Harsh Sharma, Mayuresh Patil and Craig Woolsey were supported in part by the National Science Foundation under Grant No. 1826152. Jeff Borggaard was supported in part by the National Science Foundation under grant DMS-1819110.

\bibliographystyle{vancouver}
\bibliography{main}
\end{document}